\documentclass[journal]{IEEEtran}

\newcommand{\beq}{\begin{equation}}
\newcommand{\eeq}{\end{equation}}
\newcommand{\bsq}{\begin{subequations}}
        \newcommand{\esq}{\end{subequations}}
\newcommand{\bq}{\begin{eqnarray}}
\newcommand{\eq}{\end{eqnarray}}
\newcommand{\bqn}{\begin{eqnarray*}}
        \newcommand{\eqn}{\end{eqnarray*}}
\usepackage{bbm}
\usepackage[pdftex]{graphicx}
\usepackage{enumerate}
\usepackage{enumitem} 

\usepackage{mathptmx}       
\DeclareMathAlphabet{\mathcal}{OMS}{cmsy}{m}{n}
\usepackage{helvet}         
\usepackage{courier}        
\usepackage{type1cm}        
\usepackage{makeidx}         
\usepackage{graphicx}        
\usepackage{url} 
\usepackage{multicol}        
\usepackage[bottom]{footmisc}
\usepackage{ifthen}
\usepackage{subfigure}
\usepackage[usenames,dvipsnames]{color}

\makeindex             

\usepackage{graphics} 
\usepackage{graphicx} 
\usepackage{epsfig} 
\usepackage{epstopdf}
\usepackage{mathptmx} 
\usepackage{times} 
\usepackage[cmex10]{amsmath}

\usepackage{mathtools}  
\usepackage{amssymb}  
\usepackage{amsthm}

\usepackage{amsfonts}
\usepackage{cite}
\usepackage{verbatim}
\usepackage{comment}
\usepackage{algorithmic}
\usepackage{multirow}
\usepackage{cite}
\usepackage{stfloats}
\usepackage{supertabular}
\usepackage{longtable}
\usepackage{tablefootnote}

\DeclareGraphicsExtensions{.pdf,.png,.jpg,.eps,.ps}
\renewcommand{\arraystretch}{1.2}
\usepackage{arydshln}
\usepackage{multicol}
\usepackage{colortbl}
\theoremstyle{definition}

\newtheorem{proposition}{Proposition}

\theoremstyle{definition}
\newtheorem{definition}{Definition}

\ifCLASSINFOpdf
\else
\fi

\usepackage{comment}
\usepackage{ifthen}
\newboolean{showcomments}
\setboolean{showcomments}{true}
\newcommand{\ychen}[1]{\ifthenelse{\boolean{showcomments}}
        { \textcolor{red}{YC: #1}}}
\newcommand{\tongxin}[1]{\ifthenelse{\boolean{showcomments}}
        { \textcolor{blue}{(#1)}}{}}

\usepackage{amsmath}
\allowdisplaybreaks[4]
\usepackage{algorithm}

\begin{document}
        
        %
\title{Flexibility Requirement when Tracking Renewable Power Fluctuation with Peer-to-Peer Energy Sharing}
        %
        %
        %
        
\author{Yue Chen, Wei Wei, Mingxuan Li, Laijun Chen, and Jo{\~a}o P. S. Catal{\~a}o
\thanks{This work has been submitted to the IEEE for possible publication. Copyright may be transferred without notice, after which this version may no longer be accessible.}
\thanks{Y. Chen is with the Department of Mechanical and Automation Engineering, the Chinese University of Hong Kong, Hong Kong SAR, China. (e-mail: yuechen@mae.cuhk.edu.hk)}
\thanks{W. Wei, M. Li, L. Chen are with the State Key Laboratory of Power Systems, Department of Electrical Engineering, Tsinghua University, Beijing. (e-mail: wei-wei04@mails.tsinghua.edu.cn, lmx20@mails.tsinghua.edu.cn, chenlaijun@tsinghua.edu.cn)}
\thanks{J. P. S. Catal{\~a}o is with the Faculty of Engineering of the University of Porto and INESC TEC, Porto 4200-465, Portugal (e-mail: catalao@fe.up.pt)}}

        %
        %

        \markboth{Journal of \LaTeX\ Class Files,~Vol.~XX, No.~X, Feb.~2019}%
        {Shell \MakeLowercase{\textit{et al.}}: Bare Demo of IEEEtran.cls for IEEE Journals}
        %



        \maketitle
        
        \begin{abstract}
Flexible load at the demand-side has been regarded as an effective measure to cope with volatile distributed renewable generations. To unlock the demand-side flexibility, this paper proposes a peer-to-peer energy sharing mechanism that facilitates energy exchange among users while preserving privacy. We prove the existence and partial uniqueness of the energy sharing market equilibrium and provide a centralized optimization to obtain the equilibrium. The centralized optimization is further linearized by a convex combination approach, turning into a multi-parametric linear program (MP-LP) with renewable output deviations being the parameters. The flexibility requirement of individual users is calculated based on this MP-LP. To be specific, an adaptive vertex generation algorithm is established to construct a piecewise linear estimator of the optimal total cost subject to a given error tolerance. Critical regions and optimal strategies are retrieved from the obtained approximate cost function to evaluate the flexibility requirement. The proposed algorithm does not rely on the exact characterization of optimal basis invariant sets and thus is not influenced by model degeneracy, a common difficulty faced by existing approaches. Case studies validate the theoretical results and show that the proposed method is scalable.
        \end{abstract}
        

        %
        \IEEEpeerreviewmaketitle

\section*{Nomenclature}
\addcontentsline{toc}{section}{Nomenclature}
\subsection{Indices, Sets, and Functions}
\begin{IEEEdescription}[\IEEEusemathlabelsep\IEEEsetlabelwidth{${\underline P _{mn}}$,${\overline P _{mn}}$}]
\item[$\mathcal{I}$]   Set of consumers.
\item[$\mathcal{J}$]   Set of prosumers.
\item[$l \in \mathcal{L}$]   Line $l$ in set $\mathcal{L}$.
\item[$f_k(.)$] Disutility function of user $k \in \mathcal{I} \cup \mathcal{J}$.
\end{IEEEdescription}

\subsection{Parameters}
\begin{IEEEdescription}[\IEEEusemathlabelsep\IEEEsetlabelwidth{${\underline P _{mn}}$,${\overline P _{mn}}$}]
\item[$I$] Number of consumers.
\item[$J$] Number of prosumers.
\item[$d_k$] Contract demand of user $k \in \mathcal{I} \cup \mathcal{J}$.
\item[$w_k$] Forecast renewable output of prosumer $k \in \mathcal{J}$.
\item[$\Delta w_k$] Real-time renewable output deviation of prosumer $k \in \mathcal{J}$.
\item[$\tau$] Parameter in user's objective function.
\item[$\pi_{kl}$] Line flow distribution factors.
\item[$F_l$] Power flow limit of line $l \in \mathcal{L}$.
\item[$\alpha_k,\beta_k,\zeta_k$] Parameters of the disutility function $f_k(.)$.
\item[$\underline{D}_k,\overline{D}_k$] Lower/Upper bound of demand adjustable range.
\end{IEEEdescription}

\subsection{Decision Variables}
\begin{IEEEdescription}[\IEEEusemathlabelsep\IEEEsetlabelwidth{${\underline P _{mn}}$,${\overline P _{mn}}$}]
\item[$\Delta d_k$] Demand adjustment of user $k \in \mathcal{I} \cup \mathcal{J}$.
\item[$q_k^c$] Actual sharing amount of user $k \in \mathcal{I} \cup \mathcal{J}$.
\item[$q_k$] Expected sharing amount of user $k \in \mathcal{I} \cup \mathcal{J}$.
\item[$\delta_k$] The gap between actual and expected sharing amount of user $k \in \mathcal{I} \cup \mathcal{J}$.
\item[$r_{k}^d$] Flexibility requirement of user $k \in \mathcal{I} \cup \mathcal{J}$.
\end{IEEEdescription}

\section{Introduction}
        %
        %
        %
        %
\IEEEPARstart{E}{xploiting} distributed renewable generation is an effective remedy to reduce the dependence on fossil fuel energy and achieve a sustainable society \cite{chen2020decentralized}. Meanwhile, how to tackle the volatile and intermittent energy supply caused by renewable generation has become a major concern. Existing literature concentrates on two issues: 1) how to balance the real-time power in an optimal way facing the uncertainties; 2) how to quantify the system's potential in accommodating uncertain renewable power.

\emph{For the first question}, plenty of works develop effective scheduling methods for the bulk power system that runs in a centralized manner. Typical techniques are stochastic optimization (SO) \cite{shuai2018stochastic}, robust optimization (RO) \cite{jiang2011robust}, and distributionally robust optimization (DRO) \cite{wei2015distributionally}. In stochastic optimization, the uncertain factors are modeled by their empirical distribution, and then a stochastic programming or chance-constrained problem is built. Though historical data can provide some hints for the construction of empirical distributions, the exact distribution can be hardly obtained. This inaccuracy leads to a sub-optimal scheduling strategy. Robust optimization utilizes an uncertainty set that contains all possible scenarios of the renewable output. The worst-case performance is optimized via Benders Decomposition \cite{bertsimas2012adaptive} or Column \& Constraint Generation (C\&CG) \cite{zeng2013solving} algorithms. Despite its convenience, the RO method can be too conservative since it treats all possible outputs with equal probability and neglects the fact that severe events rarely happen. Distributionally robust optimization is in-between SO and RO, where uncertainty is described by a family of inexact probability distributions restricted in an ambiguity set.

As the renewable energy generations change from large-capacity centralized units at the transmission level to small-capacity distributed units at the distribution level, various researchers seek to balance the real-time power with the help of peer-to-peer energy trading \cite{chen2018analyzing}.
Existing energy sharing mechanisms can be categorized into cooperative-game-based ones, noncooperative-game-based ones, and optimization-based ones \cite{chen2021peer}. For the first category, allocation rules are designed so that all agents have the incentive to form coalitions and together acting towards the social optimum. A coalition formation game for peer-to-peer energy trading was established with proof of several nice properties \cite{tushar2020coalition}. Network constraint \cite{azim2021coalition} and sustainable user participation \cite{tushar2019grid} were further taken into account. Though high social efficiency can be achieved, private information is needed to design the allocation rules. For the second category, two typical models are Stackelberg games and (generalized) Nash games. In Stackelberg games, the operator moves first to determine the market prices and then the users follow as price-takers. The energy sharing among prosumers was modeled as a Stackelberg game \cite{liu2017energy} and the existence of a unique and stable equilibrium was proved \cite{tushar2018peer}. A hybrid approach using stochastic optimization and Stackelberg game was developed for energy trading across the day-ahead and real-time periods \cite{liu2018energy,xu2020two}. User’s flexibility can be limited as they are price-takers under this setting. Moreover, it is hard to decide on an effective energy price especially with a large number of users, since private information may be needed and each user’s capacity is too small to be observed. In (generalized) Nash games, the users first submit their bids and then the operator clears the market. The energy trading was modeled as a two-stage stochastic Nash game \cite{li2017risk}. A generalized demand function based energy sharing mechanism was proposed \cite{chen2019energy}, and a practical bidding process as well as its convergence condition was established \cite{chen2020approaching}. The above models are for node-level sharing and how to incorporate network constraints remains to be investigated. For the third category, references \cite{cui2019peer} and \cite{cui2019efficient} established two-stage energy sharing schemes based on alternating direction method of multipliers (ADMM) algorithms. They require the help of dual variables and the economic intuition behind is hard to explain. This paper aims to develop an energy sharing mechanism based on generalized Nash game, which can protect users’ privacy, take into account network constraints, has a clear economic interpretation, and enlarge users’ flexibility by letting them be price-makers instead of price-takers. In this way, the flexibility of users is enhanced and it can facilitate the energy exchange among a wider area while network constraint is non-negligible.

\emph{For the second question}, methods to quantify the flexibility of a system are based either on regions or metrics. For regions, the do-not-exceed limit (DNEL) region \cite{zhao2014variable} and the dispatchable region \cite{wei2014dispatchable} are two well-known concepts. The DNEL region is a hypercube of uncertain parameters that will not cause infeasibility of the scheduling problem, whose size is also optimized by the problem itself. Vast literature improves the DNEL region by considering the historical data \cite{li2016data}, inexact probability distribution \cite{ma2019distributionally}, and corrective topology control \cite{korad2015enhancement}. Projecting the feasible set of the scheduling problem on the uncertainty subspace, we can get the dispatchable region. In general, the DNEL region is an inner box approximation of the dispatchable region. For metrics, the flexibility envelope \cite{nosair2015reconstructing}, the power/energy capacity \cite{dvorkin2014assessing}, etc. were proposed. A method providing both region and metric was proposed \cite{wei2020nodal}. The above works focus on centralized operation and cannot be directly applied to the energy sharing market. Reference \cite{chen2021energy} proposed the concept of absorbable region, but only region information is given and it viewed flexibility from a system level. This paper quantifies the flexibility requirement of individual users under peer-to-peer energy sharing by offering both region and metric information.

Quantifying the flexibility requirement of individual users involves solving a multi-parametric linear programming (MP-LP). In the majority of the current works, multi-parametric programming algorithms rely on the computation of critical regions where degeneracy is a key challenge. To overcome this difficulty, reference \cite{filippi2004algorithm} resorted to developing an approximate algorithm, which was then generalized to multi-parametric convex programs \cite{bemporad2006algorithm}. A prominent feature of the approximation methods is that the critical region is represented by simplices in the $p$-dimensional parameter space with $p+1$ extreme points. Such a partition may hide the impact of parameters on the optimal bases characterized by the exact critical regions which are not simplices in general. This paper proposes an alternative approximate algorithm for solving MP-LP. The basic idea is similar to \cite{filippi2004algorithm} but no simplex is needed. The main contributions are two-fold:

1) \textbf{Energy Sharing Mechanism.} A peer-to-peer energy sharing mechanism that coordinates the energy exchange among users is presented. Each user bids the quantity it would like to share aiming to minimize its disutility and match the operator’s sharing profile. The operator decides on the optimal sharing profile according to users’ bids and the physical constraints. Compared to existing works, our mechanism is easy to implement without requiring private information; can incorporate the coupling network constraints; enlarges users’ flexibility by letting them be price-makers instead of price-takers; and have several nice provable properties. Under the proposed mechanism, all users play a generalized Nash game. We prove the existence and partial uniqueness of the generalized Nash equilibrium (GNE) theoretically. The equilibrium can be retrieved from a quadratic optimization problem parameterized in the renewable output deviation, which is then turned into a multi-parametric linear program.

2) \textbf{Flexibility Characterization.} We characterize the flexibility under energy sharing based on the multi-parametric linear program derived above. To be specific, we develop an adaptive vertex generation (AVG) algorithm to partition the set of renewable output deviations into critical regions, so that in each region the demand adjustment of a user at equilibrium is shown to be a piecewise affine function in the deviation. Based on this unique feature, the flexibility requirement of individual users can be calculated. Compared to existing works, our method can provide both geometry information (critical region) and metric (flexibility requirement) about flexibility. The proposed AVG algorithm does not rely on the exact characterization of optimal basis invariant sets and thus is not influenced by model degeneracy, a common difficulty encountered by solving multi-parametric linear programs.

The rest of this paper is organized as follows. Section II presents the energy sharing mechanism and the multi-parametric linear program to calculate the equilibrium. An adaptive vertex generation algorithm is developed in Section III to solve the multi-parametric program model; based on the explicit demand adjustment policy, the flexibility requirement of individual users can be obtained. Case studies in Section IV validate the proposed models and methods. Conclusions are drawn in Section V.

\section{Peer-to-Peer Energy Sharing Model}
We consider a microgrid of $I$ consumers (indexed by $k \in \mathcal{I}=\{1,...,I\}$) and $J$ prosumers (indexed by $k \in \mathcal{J}=\{1,...,J\}$). All users $k \in \mathcal{I} \cup \mathcal{J}$ have demand while the prosumers are also equipped with renewable generators. We focus on the energy balancing at the real-time stage [C7]. At the day-ahead stage, there is renewable output predictions $w_k,\forall k \in \mathcal{I}$  based on which the users decide on their contract demands $d_k,\forall k \in \mathcal{I} \cup \mathcal{J}$. Then at the real-time stage, the actual renewable generator outputs may deviate from their predictions by $\Delta w_k,\forall k \in \mathcal{J}$. To maintain energy balancing, the users $k \in \mathcal{I} \cup \mathcal{J}$ can adjust their real-time demands by $\Delta d_k$ and exchange $q_k^c$ with each other via an energy sharing market.
The disutility caused by demand adjustment can be represented as a quadratic function $f_k(\Delta d_k):=\alpha_k (\Delta d_k)^2 + \beta_k \Delta d_k+\zeta_k$, where $\alpha_k$, $\beta_k$, and $\zeta_k$ are given parameters. This paper aims to quantify how much flexibility of individual users will be needed so that the real-time renewable generator output deviations can be balanced.

In the conventional electricity market, the penalty mechanism is implemented to encourage the participants to reveal their true production/consumption capacities and to fulfill the contract energy. This mechanism works well when the participants can accurately predict their capacities. However, considering the uncertain nature of DERs, implementing such a mechanism could lead to a huge financial loss in small prosumers and discourage them from participating in energy markets \cite{ghorani2019main}. One solution to this problem is to mitigate the overall uncertainty via local energy exchange among a group of prosumers and responsive consumers, which is what the proposed energy sharing market intends to do.

\subsection{Energy sharing mechanism}
\begin{figure}[t]
\centering
\includegraphics[width=0.7\columnwidth]{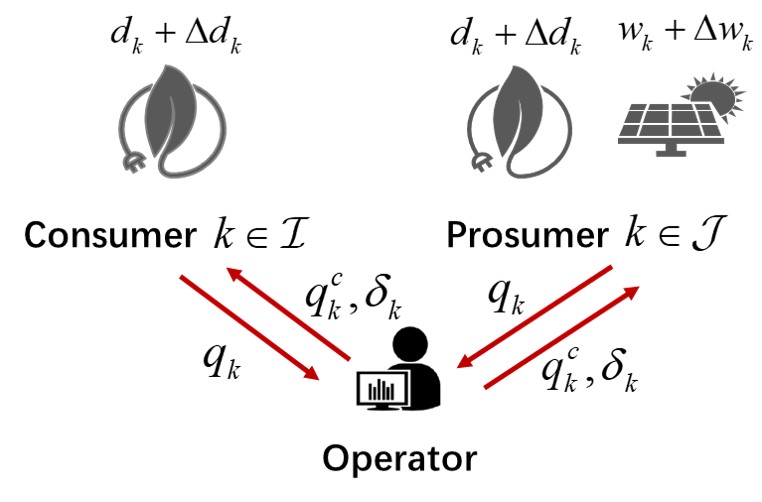}
\caption{Market structure of peer-to-peer energy sharing. Each consumer/ prosumer offers a bid $q_k$ to the operator; the operator clears the market and returns the sharing profile $q_k^c$ with the gap $\delta_k$. The users and operators are coordinated iteratively.}
\label{fig:market structure}
\end{figure}
The structure of the proposed energy sharing mechanism is shown in Fig. \ref{fig:market structure}. First, each user enters its private information into the smart meter it links to, including its disutility function $f_k(\Delta d_k)$, the contract demand $d_k$, the demand adjustable range $[\underline{D}_k,\overline{D}_k]$, and the predicted output of renewable generators $w_k$ if it is a prosumer. Then the smart meter will determine the optimal bid $q_k$ by solving \eqref{eq:consumer} or \eqref{eq:prosumer} and submit the bid to the operator. Upon receiving all the bids, the operator decides on the optimal sharing schedule to meet users' expectation as much as possible while satisfying the physical constraints by solving \eqref{eq:operator}. Then the obtained optimal schedule $q_k^c$ and unfulfilled quantity $\delta_k:=q_k^c-q_k$ are returned to the smart meter. The smart meter will adjust its bid and submit it to the operator again. This happens until the optimal schedule $q_k^c,\forall k \in \mathcal{I} \cup \mathcal{J}$ changes little between the last two iterations. Finally, the users will execute the sharing schedule $q_k^c,\forall k \in \mathcal{I} \cup \mathcal{J}$.
If $q_k^c>0$, user $k$ is a buyer; if $q_k^c<0$, user $k$ is a seller; if $q_k^c=0$, it does not take part in sharing. The key of the energy sharing mechanism design lies in how each user chooses the optimal bid $q_k$ and how the operator decides on the optimal sharing schedule $q^c_k,\forall k \in \mathcal{I} \cup \mathcal{J}$. 

For the operator, given  users' bids $q_k,\forall k \in \mathcal{I} \cup \mathcal{J}$, it solves problem \eqref{eq:operator} to set the optimal sharing schedule $q_k^c,\forall k \in \mathcal{I} \cup \mathcal{J}$:
\bsq
\label{eq:operator}
\begin{align}
    \mathop{\min}_{q_k^c,\delta_k,\forall k \in \mathcal{I} \cup \mathcal{J}} ~ & \sum \nolimits_{k \in \mathcal{I} \cup \mathcal{J}} \delta_k^2 \label{eq:operator.1}\\
    \mbox{s.t.}~ & \delta_k =q_k^c- q_k,\forall k \in \mathcal{I} \cup \mathcal{J} \label{eq:operator.2}\\
    ~ & \sum \nolimits_{k \in \mathcal{I} \cup \mathcal{J}} q_k^c = 0 \label{eq:operator.3}\\
    ~ & -F_l \le \pi_{kl}q_k^c \le F_l, \forall l \in \mathcal{L} \label{eq:operator.4}
\end{align}
\esq
The objective function \eqref{eq:operator.1} minimizes the gap between the actual and expected sharing amount.
Constraint \eqref{eq:operator.3} means the energy sold equals the energy bought. The exchanged energy is transmitted via the power network limited by constraint \eqref{eq:operator.4}.

The smart meter decides on its optimal bid $q_k$ by solving

For consumer $k \in \mathcal{I}$:
\bsq
\label{eq:consumer}
\begin{align}
    \mathop{\min}_{\Delta d_k,q_k} ~ & f_k(\Delta d_k) + \frac{\tau}{2} (q_k^c-q_k)^2 \label{eq:consumer.1}\\
    \mbox{s.t.} ~ &  q_k+\delta_k = d_k + \Delta d_k \label{eq:consumer.2}\\
    ~ & \underline{D}_k \le d_k+\Delta d_k \le \overline{D}_k \label{eq:consumer.3}
\end{align}
\esq

For prosumer $k \in \mathcal{J}$:
\bsq
\label{eq:prosumer}
\begin{align}
    \mathop{\min}_{\Delta d_k,q_k} ~ & f_k(\Delta d_k) + \frac{\tau}{2} (q_k^c-q_k)^2 \label{eq:prosumer.1}\\
    \mbox{s.t.} ~ &  w_k+\Delta w_k+q_k+\delta_k = d_k + \Delta d_k \label{eq:prosumer.2}\\
    ~ & \underline{D}_k \le d_k+\Delta d_k \le \overline{D}_k \label{eq:prosumer.3}
\end{align}
\esq
The objective of each user $k \in \mathcal{I} \cup \mathcal{J}$ consists of two parts: its disutility $f_k(\Delta d_k)$ and the gap between its bid $q_k$ and operator's schedule $q_k^c$. Parameter $\tau$ indicates the trade-off between the above two parts and help evaluate the gap in monetary terms.
 Constraints include the power balance conditions for consumers \eqref{eq:consumer.2} and for prosumers \eqref{eq:prosumer.2}, and the demand adjustable range limits \eqref{eq:consumer.3} and \eqref{eq:prosumer.3}.
 
The proposed energy sharing mechanism well fits the practice in that the users do not need to know the network constraint \eqref{eq:operator.4} available only to the operator, and the operator does not require local constraint \eqref{eq:consumer.3} and \eqref{eq:prosumer.3} private to each user. Compared to centralized dispatch, our mechanism has two advantages: 1) \textbf{Regarding privacy.} Centralized dispatch needs to gather information from all users to make the central decision. However, the users may be unwilling to provide this sort of private information. In the proposed peer-to-peer energy sharing, each user only needs to submit a bid to the operator with its private data input only to its own smart meter, so privacy can be preserved.
2) \textbf{Less communication burden.} Under peer-to-peer energy sharing, the only information exchanged will be the bid, which is a scalar. However, under centralized dispatch, the operator needs to collect information such as cost coefficients, demands, and demand adjustable boundaries, from all users, which incurs a high communication burden overhead.  Meanwhile, this mechanism can achieve a social optimal outcome (proved in Proposition \ref{prop-1} later).

\textbf{Remark 1:} Our proposed energy sharing mechanism focuses on the real-time market \cite{pei2016optimal}, which will be conducted hourly. In real-time, when the uncertainty is realized, we resort to a mechanism similar to the economic dispatch \cite{wei2014robust} (happens every several mins or an hour) to balance the system. Specially, instead of balancing the real-time power via centralized dispatch, we develop an energy sharing mechanism that allows the users to adjust their elastic demands and exchange energy with each other. AGC that happens at a smaller time resolution and usually at the transmission level is not considered in this paper as most of the economic dispatch works \cite{wei2014robust,wei2015distributionally}.

\textbf{Remark 2:} In this paper, deciding on the optimal bids and finally the sharing quantities involve solving optimization problems \eqref{eq:consumer} or \eqref{eq:prosumer} and two-way communication with the operator. With the current technology, advanced smart meters have certain data processing abilities and enable two-way communication between the meter and the central system \cite{nimbargi2016review}. Therefore, it is reasonable to assume that with the assist of a smart meter, the user would be able to decide on the optimal bids. When taking into account various kinds of appliances and facilities, the decision-making may become more complex. In that case, we may need a more sophisticated energy management system to handle all the devices.

\subsection{Market equilibrium}
\begin{table*}[h]
        \renewcommand{\arraystretch}{1.3}
        \renewcommand{\tabcolsep}{1em}
        \centering
        \caption{Element of the energy sharing game}
        \label{tab:game}
        \begin{tabular}{ccc}
                \hline 
               Players $\mathcal{K}$ & Actions $\mathcal{X}$ & Cost functions $\Pi$\\
               \hline
               consumer $k=1,\cdots, I$ & $\mathcal{X}_k(q_k^c,\delta_k):=\{(\Delta d_k,q_k)~|~ \mbox{\eqref{eq:consumer.2} and \eqref{eq:consumer.3} are met.} \}$ & $\Pi_k(q_k^c,\delta_k):=f_k(\Delta d_k) + \tau/2 (q_k^c-q_k)^2$\\
               prosumer $k=(I+1),\cdots, (I+J)$ & $\mathcal{X}_k(q_k^c,\delta_k):=\{(\Delta d_k,q_k)~|~ \mbox{\eqref{eq:prosumer.2} and \eqref{eq:prosumer.3} are met.} \}$ & $\Pi_k(q_k^c,\delta_k):=f_k(\Delta d_k) + \tau/2 (q_k^c-q_k)^2$\\
               operator $k=I+J+1$ & $\mathcal{X}_k(\Delta d,q):=\{(q^c,\delta)~|~ \mbox{\eqref{eq:operator.2} and \eqref{eq:operator.3} are met.} \}$ &  $\Pi_{I+J+1}:= \sum \nolimits_{k \in \mathcal{I} \cup \mathcal{J}} \delta_k^2$\\
                \hline
        \end{tabular}
\end{table*}

The above peer-to-peer energy sharing market can be regarded as a generalized Nash game \cite{harker1991generalized} with its elements summarized in TABLE \ref{tab:game}. For simplicity, we use $\mathcal{G}=\{\mathcal{K},\mathcal{X},\Pi\}$ to denote the sharing game in an abstract form, where $\mathcal{K}$ is the set of players, $\mathcal{X}$ the action set, and $\Pi$ the collection of cost functions. It's worth noting that, the action sets of users depend on the strategies of the operator (the $q_k^c, \delta_k,\forall k \in \mathcal{I} \cup \mathcal{J}$); the action sets of the operator depends on the strategies of the users (the $\Delta d_k, q_k,\forall k \in \mathcal{I} \cup \mathcal{J}$). Therefore, it constitutes a generalized Nash game, whose equilibrium is harder to analyze than the standard Nash game due to the complex coupling \cite{facchinei2010generalized}. The equilibrium of the generalized Nash game, from which every player has no incentive to deviate, is defined below.
\begin{definition}(Generalized Nash Equilibrium) A profile $(\Delta d^*, q^*, q^{c*}, \delta^*) \in \mathcal{X}$ is a \emph{generalized Nash equilibrium} (GNE) of the energy sharing game $\mathcal{G}=\{\mathcal{K},\mathcal{X},\Pi\}$ if $\forall k=1,...,(I+J)$
\begin{align}
\label{eq:game1}
    (\Delta d_k^*,q_k^*) = ~ &  \mbox{argmin}_{\Delta d_k,q_k} \Pi_k(q_k^{c*},\delta_k^*) \nonumber\\
    ~ & \mbox{s.t.} ~~ (\Delta d_k,q_k) \in \mathcal{X}_k(q_k^{c*},\delta_k^*)
\end{align}
and
\begin{align}
\label{eq:game2}
    (q^{c*},\delta^*) = ~ &  \mbox{argmin}_{q^c,\delta}  \Pi_{I+J+1} \nonumber\\
    ~ & \mbox{s.t.}~~ (q^c,\delta) \in \mathcal{X}_{I+J+1} (\Delta d^*,q^*)
\end{align}
\end{definition}
 
 In the following proposition, we prove the existence of the equilibrium and offer a centralized optimization problem to compute it, which casts down to multi-parametric program parameterized in the renewable output deviations. Denote $\mathcal{Q}:=\{q^c~|~ \mbox{\eqref{eq:operator.3} and \eqref{eq:operator.4} are met}\}$, $\mathcal{D}_k:=\{\Delta d_k~|~ \underline{D}_k \le d_k+\Delta d_k \le \overline{D}_k\}$ for all $k \in \mathcal{I} \cup \mathcal{J}$.

\begin{proposition}
\label{prop-1}
The energy sharing game $\mathcal{G}=\{\mathcal{K},\mathcal{X},\Pi\}$ has at least one GNE if \eqref{eq:central} is feasible, and $(\Delta d^*,q^*,q^{c*},\delta^*)$ is an GNE if and only if $\Delta d^*$ is the unique optimal solution of 
\bsq
\label{eq:central}
\begin{align}
    \mathop{\min}_{\Delta d_k,\forall k \in \mathcal{I} \cup \mathcal{J}} ~ & \sum \nolimits_{k \in \mathcal{I} \cup \mathcal{J}} f_k(\Delta d_k) \label{eq:central.1} \\
    \mbox{s.t.}~& q_k^{c} = \left\{ \begin{aligned}
    & d_k+\Delta d_k,\forall k \in \mathcal{I} \\
    & d_k+\Delta d_k-w_k-\Delta w_k ,\forall k \in \mathcal{J}
    \end{aligned}\right.:\eta_k \label{eq:central.2}\\
    ~ & \Delta d_k \in \mathcal{D}_k,\forall k \in \mathcal{I} \cup \mathcal{J}\label{eq:central.3}\\
    ~ & q^c \in \mathcal{Q} \label{eq:central.4}
\end{align}
\esq
with $\delta^*_k=-\hat \eta_k/\tau$ for all $k \in \mathcal{I} \cup \mathcal{J}$, where $\hat \eta$ is the value of dual variable at optimum, and $q_k^{c*}=q_k^*+\delta^*_k$ with
\begin{align}
q_k^*=\left\{
\begin{aligned}
& d_k+\Delta d_k^* -\delta^*_k~,~\forall k \in \mathcal{I} \\
& d_k+\Delta d_k^*-w_k-\Delta w_k  -\delta^*_k~,~\forall k \in \mathcal{J}  \\
\end{aligned}
\right.
\end{align}
\end{proposition}

The proof of Proposition \ref{prop-1} can be found in Appendix \ref{apen-1}. It offers a more convenient way to calculate and analyze the energy sharing game equilibrium by solving \eqref{eq:central}. As problem \eqref{eq:central} minimizes the total disutility of all users, $\Delta d^*$ is social optimal. Take the real-time renewable output deviation $\Delta w$ as parameters, problem \eqref{eq:central} can be regarded as a quadratic multi-parametric optimization problem. 

\subsection{Linearization and compact form}

For the sake of analysis, the objective function is linearized via a convex combination approach \cite{wu2011tighter}. Take a univariate convex function $z=g(\xi)$, $\xi_l \le \xi \le \xi_u$ for example. The function is evaluated at some discrete points $\xi_1,\cdots,\xi_K$, and $z_1=g(\xi_1),\cdots,z_K=g(\xi_K)$. By introducing variables $\sigma_1,\cdots,\sigma_K \ge 0$ and $\sum_{k=1}^K \sigma_k =1$, $\xi$ and $g(\xi)$ in the optimization problem can be replaced with linear functions $\sum_{k=1}^K \sigma_k \xi_k$ and $\sum_{k=1}^K \sigma_k z_k$ in $\sigma$, respectively. Applying this technique to \eqref{eq:central.1}, problem \eqref{eq:central} turns into a multi-parametric linear program (MP-LP) with the compact form as:
\begin{equation}
\begin{aligned}
v(\theta) = \mathop{\min}_{x} ~~ & c^\top x \\
\mbox{s.t.} ~~ & A x  \le t + B \theta  
\end{aligned}
\label{eq:mp-LP-Matrix}
\end{equation}
where $x \in \mathbb R^n$ and $\theta \in \mathbb R^p$ are vectors of decision variables and parameters, respectively; $A \in \mathbb R^{m \times n}$, $t \in \mathbb R^m$, $B \in \mathbb R^{m \times p}$ and $c \in \mathbb R^n$ are input data. The decision variables $x$ is bounded due to the lower and upper bound constraints. Denote $v(\theta)$ and $x^*(\theta)$ as the optimal value and optimal solution associated with given parameter $\theta$, respectively.

Based on the above MP-LP \eqref{eq:mp-LP-Matrix}, in the following, we develop an adaptive vertex generation (AVG) algorithm to approximate the optimal cost and the demand adjustment at equilibrium under different renewable output deviations and quantify the flexibility requirements of individual users.

\section{Flexibility Characterization}
\begin{figure}[t]
\centering
\includegraphics[width=0.45\columnwidth]{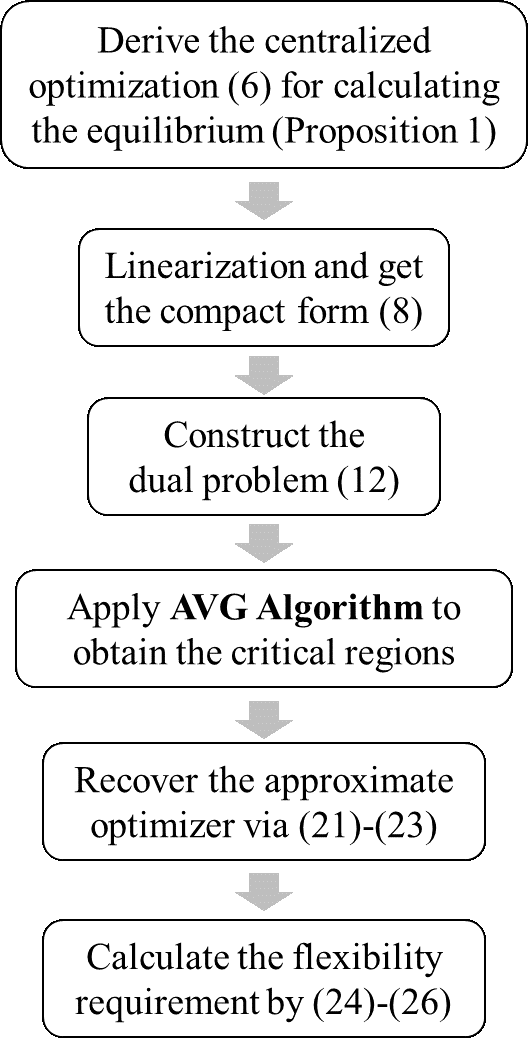}
\caption{A linear program \eqref{eq:mp-LP-Matrix} for computing the market equilibrium is derived based on Proposition 1 and linearization techniques. Then AVG algorithm is applied to calculate the flexibility requirement.}
\label{fig:framework}
\end{figure}

In this section, we first define the flexibility requirement of individual users under energy sharing; then, an adaptive vertex generation algorithm is developed to construct a piecewise linear estimator of the optimal total cost and the demand adjustment at equilibrium subject to a given error tolerance, based on which we calculate the flexibility requirement. The processing framework of the algorithms is shown in Fig. \ref{fig:framework}.

\subsection{Definition of flexibility requirement}
Let $\bar \Theta$ be the maximum set of $\theta$ such that problem (\ref{eq:mp-LP-Matrix}) has at least one feasible solution $x$, known as the dispatchable region \cite{wei2014dispatchable}. We assume $\theta \in \Theta$, a bounded polyhedral  subset of $\bar \Theta$. The flexibility requirement of each user is defined below.
\begin{definition} (Flexibility Requirement) Suppose $\Delta d_k^*(\theta),\forall k \in \mathcal{I} \cup \mathcal{J}$ are users' demand adjustment at equilibrium under renewable deviation $\theta$. The flexibility requirement of user $k \in \mathcal{I} \cup \mathcal{J}$ is defined as a region $[\underline{r}_k^d,\overline{r}_k^d]$ with
 \begin{align}
 \label{eq:definition-1}
     \underline{r}_k^d = \min \{\Delta d_k^*(\theta), \forall \theta \in \Theta\}
 \end{align}
 and
  \begin{align}
  \label{eq:definition-2}
     \overline{r}_k^d = \max \{\Delta d_k^*(\theta), \forall \theta \in \Theta\}
 \end{align}
\end{definition}

The flexibility in this paper is a kind of \emph{reserve capacity} \cite{wei2014robust}. It defines the backup adjustable capacity the demand can provide in the occurrence of renewable generation output deviations due to uncertainty. Since $\Delta d^*(\theta)$ involves solving problem \eqref{eq:central}, getting $[\underline{r}_k^d, \overline{r}_k^d]$ directly by definition can be difficult. In this paper, an MP-LP based approach is proposed to obtain the flexibility requirement. To be specific, we divide $\Theta$ into non-overlapping critical regions CR$_1$, $\cdots$, CR$_N$ such that in each region CR$_i$, the optimal value $v(\theta)$ and optimal solution $x^*(\theta)$ are linear functions in the parameter vector $\theta$. Then \eqref{eq:definition-1}-\eqref{eq:definition-2} is equivalent to finding the maximum/minimum point of a piecewise linear function, which is easier to solve.

\subsection{Adaptive vertex generation algorithm}
The majority of current MP-LP related works identifies the critical regions by leveraging the graph of optimal bases based on which the expression of $x^\star(\theta)$ and $v(\theta)$ are obtained \cite{wei2020nodal}. However, degeneracy that indicates the existence of multiple primal or dual optimal solutions may cause difficulties in building the critical regions \cite{borrelli2003geometric}. In this paper, we develop an adaptive vertex generation algorithm to approximate $v(\theta)$ via dual variables without the information of critical regions first, and then retrieve the critical region and $x^*(\theta)$ accordingly. In this way, the impact of degeneracy can be avoided. The details of the algorithm are as follows:

\emph{1) Lower bound of $v(\theta)$ from dual problem}

Suppose the optimal value function $v(\theta)$ can be represented by a piecewise linear function 
\begin{equation}
v(\theta) = \begin{cases}
m_1 + n_1^\top \theta, &  \theta \in \mbox{CR}_1 \\
m_2 + n_2^\top \theta, &  \theta \in \mbox{CR}_2 \\
\qquad \vdots  \\
m_N + n_N^\top \theta, &  \theta \in \mbox{CR}_N \\
\end{cases}
\label{eq:OVF-Deter}
\end{equation}
where $m_1,\cdots,m_N$ are constant scalars; $n_1,\cdots,n_N$ are constant vectors.
For any fixed $\theta$, the dual problem of (\ref{eq:mp-LP-Matrix}) is
\begin{equation}
\begin{aligned}
v(\theta) = \mathop{\max}_{\gamma}~~& \gamma^\top ( t + B \theta) \\
\mbox{s.t.} ~~ & A^\top \gamma = c,~ \gamma \le 0 
\end{aligned}
\label{eq:Dual-LP}
\end{equation}
where $\gamma$ is the dual variable, and its feasible region is denoted by
$\Gamma = \left\{ \gamma \middle| A^\top \gamma = c, \gamma \le 0 \right\}$. 
According to strong duality, the optimal value of (\ref{eq:Dual-LP}) is bounded and equals to $v(\theta)$ for any given $\theta$. Since $v(\theta)$ is finite, the dual optimal solution $\gamma$ must belong to the set of extreme points of $\Gamma$ [vert($\Gamma$) for short], although $\Gamma$ may be unbounded and contain extreme rays. In this regard, by enumerating the vertices of $\Gamma$, $v(\theta)$ can be equivalently expressed as 
\begin{equation}
\begin{gathered}
v(\theta) = \max_i \left\{ \gamma^{\top}_i t+ \gamma^{\top}_i B \theta \right\}\\
\{\gamma_1,\gamma_2,\cdots\} \in \mbox{vert}(\Gamma)
\end{gathered}  \label{eq:OVF-Dual}
\end{equation}
Comparing (\ref{eq:OVF-Dual}) and (\ref{eq:OVF-Deter}), coefficients $m_i$ and $n_i$ and dual variable $\gamma_i$ have the following relation 
\begin{equation}
m_i = \gamma^{\top}_i t,~ n_i = \gamma^{\top}_i B  
\label{eq:PWL-Coefficients}
\end{equation} 
Nevertheless, vertex enumeration is an exhaustive task, and only a small fraction of vertices correspond to valid pieces in (\ref{eq:OVF-Deter}) and (\ref{eq:OVF-Dual}), but we do not know which $\mbox{one in vert}(\Gamma)$ appears in  (\ref{eq:OVF-Dual}) in advance. Let $\underline \Gamma$ be a subset of $\mbox{vert}(\Gamma),$ and
\begin{equation}
\underline v(\theta) = \max_i \left\{ \gamma^{\top}_i t+ \gamma^{\top}_i B \theta \right\}, \{\gamma_1,\gamma_2,\cdots\} \in \underline \Gamma
\label{eq:OVF-Dual-L}
\end{equation}
Hence, $\underline v(\theta)$ is an underestimator of $v(\theta)$. The approximation quality depends on the choice of $\underline \Gamma$. In the following, an adaptive vertex generation algorithm is developed to construct $\underline \Gamma$ so that the gap between $\underline v(\theta)$ and $v(\theta)$ can be small.

\emph{2) Recover the critical regions}

Suppose (\ref{eq:OVF-Dual-L}) is a non-redundant representation. Intuitively, if the value of any piece $\gamma^{\top}_i t+ \gamma^{\top}_i B \theta$ can reach maximum for some $\theta \in \Theta$, (\ref{eq:OVF-Dual-L}) is non-redundant. The exact definition and redundancy elimination methods are in Appendix \ref{apen-2}. Then the critical region associated with each piece can be retrieved from (\ref{eq:OVF-Dual-L}).   
Recall the original piecewise linear form (\ref{eq:OVF-Deter}) and the point-wise maximum form (\ref{eq:OVF-Dual-L}), the value of $m_i + n_i^\top \theta$ achieves maximum in CR$_i$, i.e.: 
\begin{equation}
\label{eq:CR-OVF}
\mbox{CR}_i = \left\{ \theta \middle| m_i + n_i^\top \theta \ge m_{[-i]} + n_{[-i]}^\top \theta \right\}
\end{equation}
where $[-i]=\{1,\cdots, i-1, i+1,\cdots,n\}$ stands for the set of complementary indexes of $i$, so CR$_{i}$ is described by $n-1$ inequalities in (\ref{eq:CR-OVF}). It can be represented in a matrix form \eqref{eq:Polyhedral-CR}, where $H \in \mathbb R^{(n-1)\times p}$, $h \in \mathbb R^p$.
\begin{equation}
\mbox{CR}_i = \left\{ \theta \middle| H \theta \le h \right\}
\label{eq:Polyhedral-CR}
\end{equation}
For notation conciseness, index $i$ of critical region for matrix $H$ and vector $h$ are omitted; $H_j \theta \le h_j$, $j=1,\cdots,(n-1)$ represent the individual constraints in CR$_i$.

\emph{3) Overall algorithm}

Suppose we already have an approximation (\ref{eq:OVF-Dual-L}) for the optimal value function, as well as critical regions retrieved from (\ref{eq:CR-OVF}), then we have to investigate the quality of this approximation in each CR$_i$. The task comes down to evaluate the maximal approximation error through
\begin{equation}
\label{eq:Appr-Err-1}
\begin{aligned}
\varepsilon_i = \max ~ & v(\theta)-\underline v(\theta)  \\
\mbox{s.t.}~ &\theta \in \mbox{CR}_i 
\end{aligned}
\end{equation}
where $v(\theta)$ is the true optimal value function of problem (\ref{eq:mp-LP-Matrix}), which is unknown. $\varepsilon_i$ must be nonnegative because $\underline v(\theta)$ is an underestimator. First, due to strong duality, we have
\begin{equation}
\label{eq:proof-1}
\begin{aligned}
& \qquad v(\theta)-\underline v(\theta) \\
& = \max_{\theta \in \mbox{CR}_i} \left\{\max_{\gamma \in \Gamma} \gamma^\top ( t + B \theta) - (\gamma^{\top}_i t+ \gamma^{\top}_i B \theta) \right\}
\end{aligned}
\end{equation}

Problem (\ref{eq:proof-1}) is nonconvex due to the product term $\gamma^\top B\theta$ in the objective function. Since the parameter space has a low dimension, we can solve (\ref{eq:proof-1}) by numerating the vertex of $\mbox{CR}_i=\{\theta|H\theta \le h\}$, which is
\begin{equation}
\label{eq:proof-2}
\begin{aligned}
\epsilon_i =\max_i \max_{\gamma}~ & (\gamma -\gamma_i)^\top (t+B \theta_i ) \\
\mbox{s.t.}~ &  \{\theta_1,\theta_2,...\} \in \mbox{vert(CR}_i)\\ 
& A^\top \gamma = c,~ \gamma \le 0
\end{aligned}
\end{equation}

Another approach to solve \eqref{eq:proof-1} is the big-M method \cite{fortuny1981representation}, but our approach is more robust since it avoids the difficulty of choosing a proper $M$. The proposed adaptive vertex generation algorithm is summarized in Algorithm 1. The output is a $\varepsilon$-optimal underestimator of the true optimal value function.

\begin{algorithm}
\normalsize
\label{Ag:Path-Generation}
\caption{{\bf : Adaptive Vertex Generation Algorithm}}
\begin{algorithmic}[1]
\STATE Initiation:  Error tolerance $\varepsilon>0$; parameter set $\Theta$; initial sampled parameters  $\theta_1,\cdots,\theta_n$; 
 
\STATE Underestimation

{\bf For} $i=1:n$ 

\qquad Solve LP (\ref{eq:Dual-LP}) with $\theta_i$ and the solution is $\gamma^*_i$.

\qquad Update $\underline{\Gamma} \leftarrow \underline{\Gamma} \cup \gamma^*_i$ if $\gamma^*_i \notin \underline{\Gamma}$.

{\bf end}
 
Construct $\underline v(\theta)$ by (\ref{eq:OVF-Dual-L}); retrieve CR$_1$, $\cdots$, CR$_{n}$ by \eqref{eq:CR-OVF}.

\STATE Check approximation error

{\bf For} each critical region $i$ 

\qquad Solve LP (\ref{eq:proof-2}) with CR$_i$. Record $\varepsilon^*_i$, $\gamma^*$.

\qquad If $\varepsilon^*_i > \varepsilon$, update $\underline \Gamma \leftarrow \underline \Gamma \cup \gamma^*$.

{\bf end}

If $\varepsilon^*_i \le \varepsilon$, $\forall i$, terminate.

\STATE Update  $\underline v(\theta)$ by (\ref{eq:OVF-Dual-L}) with the current $\underline \Gamma$, and remove redundant pieces (see Appendix \ref{apen-2}). 

\STATE Retrieve critical regions from the current   $\underline v(\theta)$; remove redundant constraints (Appendix \ref{apen-2}), and go step 3.
\end{algorithmic}
\end{algorithm}

\textbf{Remark:} Because the set $\mbox{vert}(\Gamma)$ has finite elements, Algorithm 1 must terminate in a finite number of steps. Nonetheless, as only a small fraction of elements in $\mbox{vert}(\Gamma)$ corresponds to the pieces in the optimal value function, Algorithm 1 is efficient in practice. In contrast to the existing approaches where critical regions are obtained from the analysis of optimal bases invariancy,  in Algorithm 1, critical regions are retrieved from the optimal value function. Therefore, model degeneracy does not affected the proposed algorithm.

\subsection{Calculation of flexibility requirement}
To obtain the flexibility requirement of individual users, we need to recover an approximate optimizer $x(\theta)$ as a function of $\theta$. In CR$_i$, the optimal value is $v(\theta) = \gamma^{\top}_i t+ \gamma^{\top}_i B \theta$, the complementarity and slackness condition of (\ref{eq:mp-LP-Matrix}) implies
\begin{equation}
\gamma^*_{ij} < 0 \to (Ax - B\theta-t)_j = 0
\end{equation}

Let $A^\prime_i/B^\prime_i/t^\prime_i$ denote the rows of $A/B/t$ corresponding to the indexes of nonzero elements in $\gamma_i^*$, we obtain a set of linear equations:
\begin{equation}
\label{eq:re-group}
A_i^\prime x = t_i^\prime + B_i^\prime\theta,~  \theta \in \mbox{CR}_i
\end{equation}
If $A_i^\prime$ is invertible, the optimizer is obtained by   
\begin{equation}
\label{eq:optimizer-1}
x(\theta) = (A_i^\prime)^{-1} (t_i^\prime + B_i^\prime\theta),~  \theta \in \mbox{CR}_i
\end{equation}
Otherwise, we can calculate the flexibility requirement $[\underline{r}_{ik}^d,~\overline{r}_{ik}^d]$ by the following optimization problems. Without loss of generality, suppose the first $I+J$ terms in $x$ correspond to the decision variables $\Delta d_k,\forall k=1,...,(I+J)$.
\bsq
\begin{align}
    \underline{r}_{ik}^d = \mathop{\min}_{\theta \in \mbox{CR}_i} ~ & x_k \\
    \mbox{s.t.} ~ & A_i^\prime x = t_i^\prime + B_i^\prime \theta \\
    ~ & \underline{D}_k \le d_k+x_k \le \overline{D}_k, \forall k=1,...,(I+J)
\end{align}
\esq
and
\bsq
\begin{align}
    \overline{r}_{ik}^d = \mathop{\max}_{\theta \in \mbox{CR}_i} ~ & x_k \\
    \mbox{s.t.} ~ & A_i^\prime x = t_i^\prime + B_i^\prime \theta \\
    ~ & \underline{D}_k \le d_k+x_k \le \overline{D}_k, \forall k=1,...,(I+J)
\end{align}
\esq

Therefore, the flexibility requirement of user $k$ over the entire parameter set $\theta$ is $[\underline{r}_k^d,\overline{r}_k^d]$, where
\begin{align}
    \underline{r}_k^d = \mathop{\min}_{i} \{r_{ik}^d\},~ \overline{r}_k^d = \mathop{\max}_{i} \{r_{ik}^d\}
\end{align}

In general, the obtained critical regions CR$_i,\forall i$ can provide geometry information of how uncertainty be handled under energy sharing, while the flexibility requirement $[\underline{r}_k^d, \overline{r}_k^d],\forall k$ offers a quantitative metric.

\section{Case Studies}
We first use a 5-bus system for illustration; and then a 69-bus microgrid system is tested to show the scalability of the proposed model and algorithm. Data is available in \cite{Data}.

\begin{figure}[t]
\centering
\includegraphics[width=0.8\columnwidth]{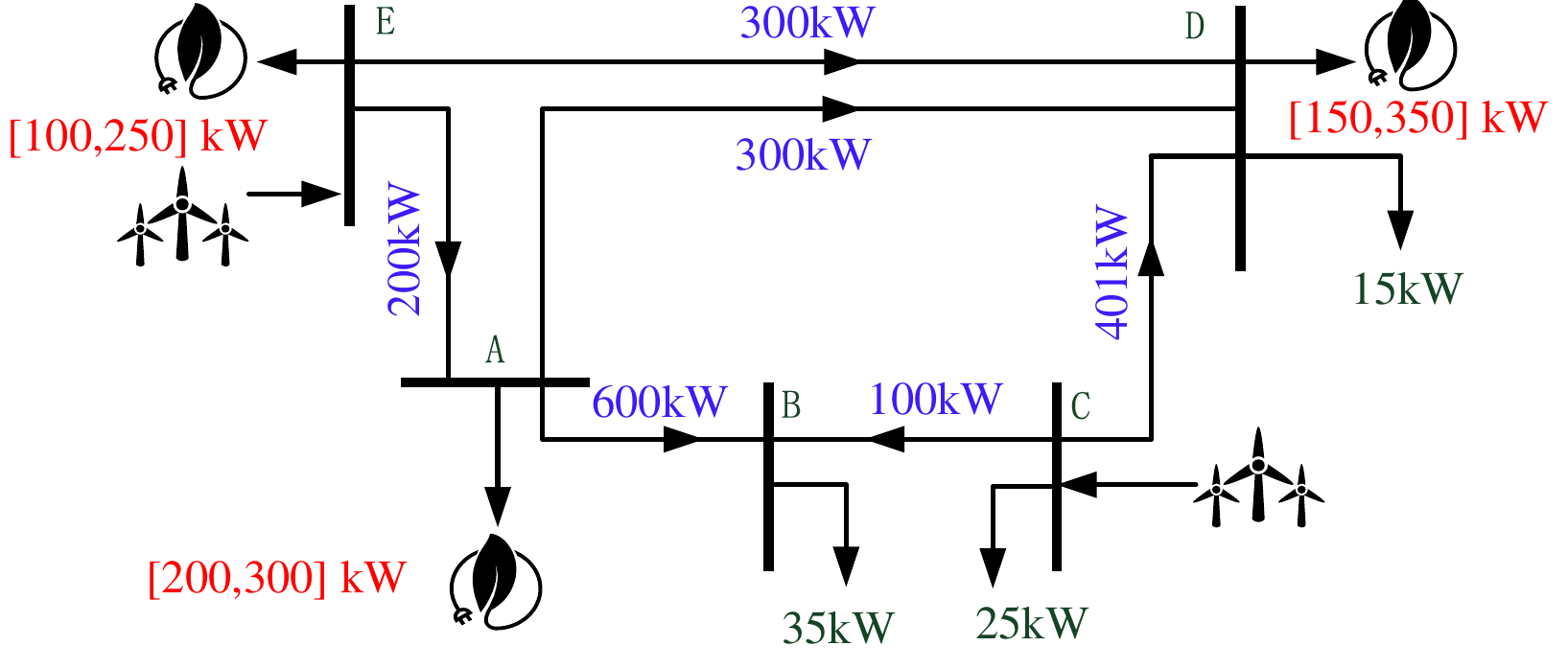}
\caption{Topology and parameters of the 5-bus system. The adjustable ranges of elastic demands are in red, the line flow limits are in blue, and the fixed demands are in dark green.}
\label{fig:PJM}
\end{figure}

\subsection{5-bus system}
\begin{table}[t]
        \renewcommand{\arraystretch}{1.3}
        \renewcommand{\tabcolsep}{1em}
        \centering
        \caption{System Parameter of the 5-bus System}
        \label{5bussystem}
        \begin{tabular}{cccc}
                \hline 
               From & To & $X$ & $F_l$\\
               \hline
               A & B & 0.0281 & 600\\
               A & D & 0.0304 & 300\\
               A & E & 0.0064 & 200\\
               B & C & 0.0108 & 100\\
               C & D & 0.0297 & 401\\
               D & E & 0.0297 & 300\\
                \hline
        \end{tabular}
\end{table}

\begin{table}[t]
        \renewcommand{\arraystretch}{1.3}
        \renewcommand{\tabcolsep}{0.7em}
        \centering
        \caption{Demand Data of the 5-bus System\tablefootnote{Wind farms locate at bus 3, 5 with $w_1=220$ kW, $w_2=450$ kW.}}
        \label{tab:5busdemand}
        \begin{tabular}{cccccc}
                \hline 
             Node & $d_k$ (kW)  & Bound & $\alpha_k$ (\$/$\mbox{kW}^2$) & $\beta_k$ (\$/kW) & $\zeta_k$ (\$)\\
               \hline
               A & 230 & [200, 300] & 0.003 & 1.80 & 255.30\\
               B & 35 &	- &	- & - & -\\
               C & 25 & - & - & - & -\\
               D & 185 \tablefootnote{170 kW elastic demand and 15 kW inelastic demand.} & [150, 350] & 0.006 & 2.76 & 295.80\\
               E & 200 & [100, 250] & 0.005 & 2.56 & 312.00\\
                \hline
        \end{tabular}
\end{table}

\begin{figure}[t]
\centering
\includegraphics[width=0.7\columnwidth]{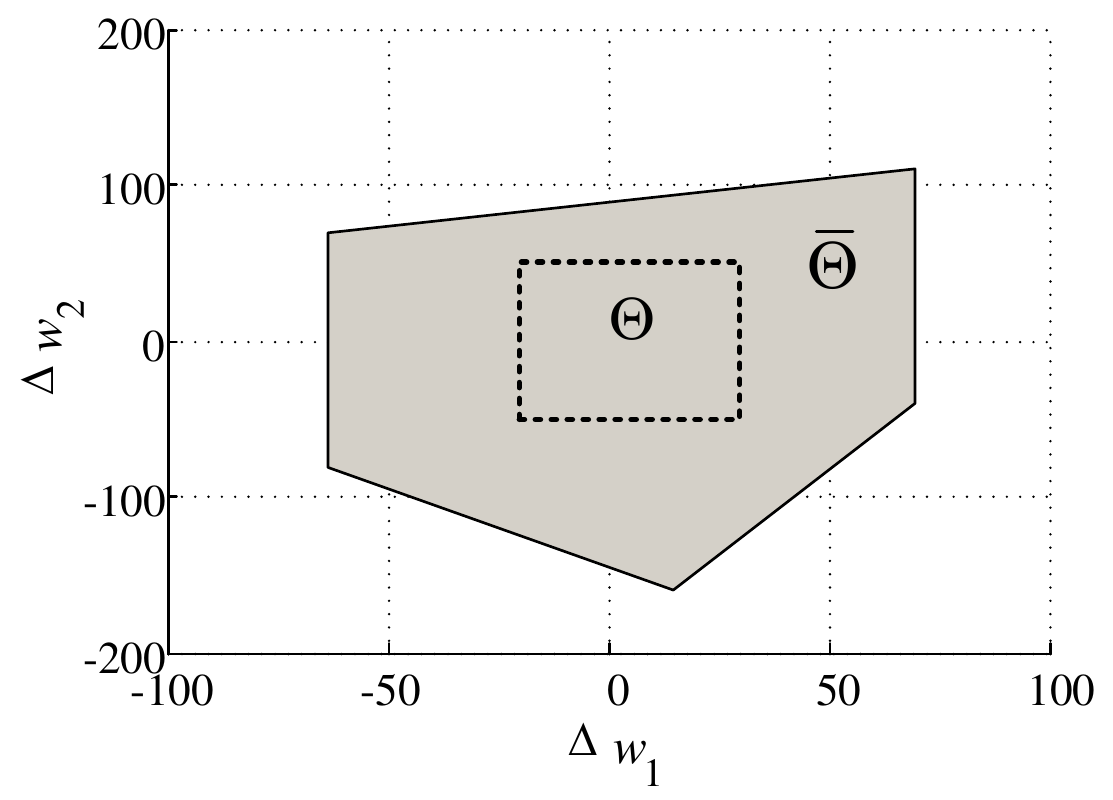}
\caption{Dispatchable region $\bar \Theta$ (maximum set of $\theta$ such that problem \eqref{eq:mp-LP-Matrix} has at least one feasible solution) and parameter range $\Theta$ (we assume $\theta \in \Theta$) of the 5-bus system.}
\label{fig:dispatchable}
\end{figure}

\begin{figure}[t]
\centering
\includegraphics[width=0.7\columnwidth]{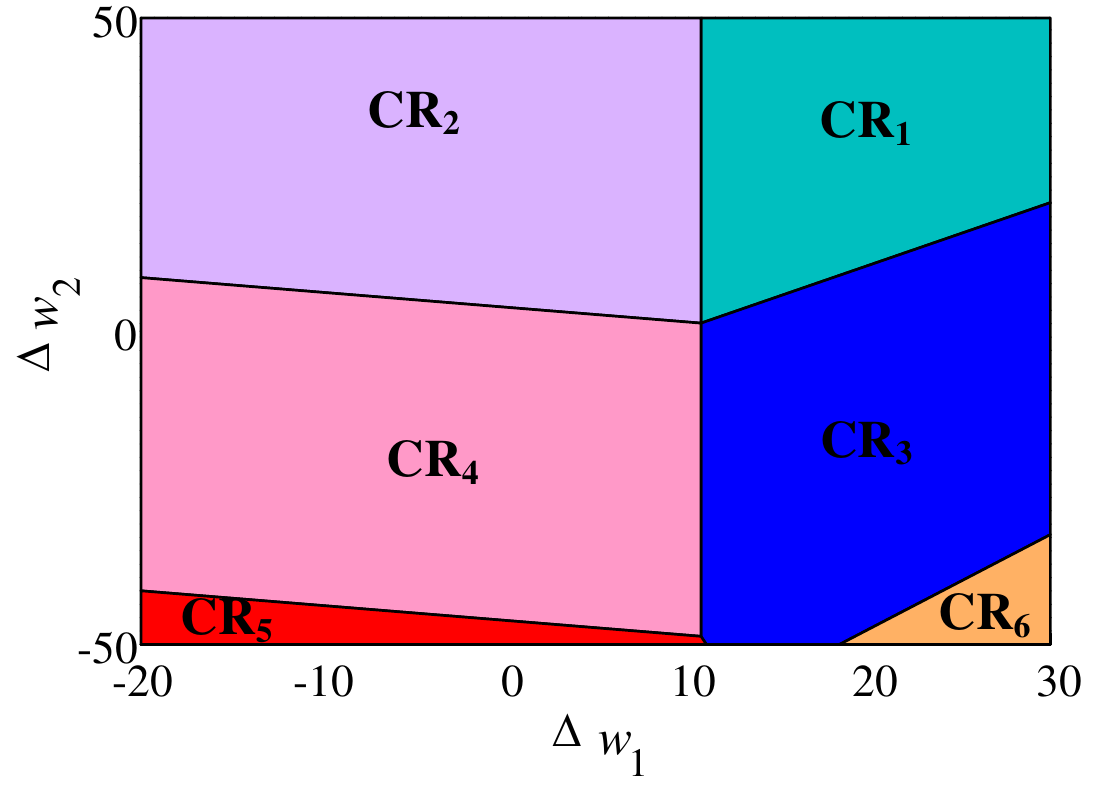}
\caption{Critical Regions of the 5-bus system. The parameter range is divided into six critical regions and within each region, the optimal cost is a linear function in $\theta$ given in \eqref{eq:approximate-obj}.}
\label{fig:critical}
\end{figure}

The topology of the 5-bus system is given in Fig. \ref{fig:PJM} with other parameters in TABLE \ref{5bussystem}-\ref{tab:5busdemand}. There are three users with elastic demands and two renewable generators at nodes C and E, respectively, so $\theta=[\Delta w_1, \Delta w_2]^\top$. The maximum set of parameters $\bar \Theta$ is the dispatchable region \cite{wei2014dispatchable}, the grey region in Fig. \ref{fig:dispatchable}. The $\theta$ can vary within the rectangle subregion $\Theta$. Applying Algorithm 1, we can obtain six critical regions as in Fig. \ref{fig:critical}. Note that the quadratic objective function \eqref{eq:central} is linearized by 5 segments. The approximate optimal cost function $\underline v(\theta)$ in each critical region is as follows.
\begin{equation}
\label{eq:approximate-obj}
v(\theta) = \begin{cases}
814.92 + 4.00\Delta w_1+2.81 \Delta w_2, &  \theta \in \mbox{CR}_1 \\
835.39 + 2.13\Delta w_1 + 2.81\Delta w_2, &  \theta \in \mbox{CR}_2 \\
810.04 + 4.50\Delta w_1 + 2.31\Delta w_2, &  \theta \in \mbox{CR}_3 \\
837.26 + 2.01\Delta w_1 + 2.31\Delta w_2, &  \theta \in \mbox{CR}_4 \\
819.08 + 1.92\Delta w_1 + 1.92\Delta w_2, &  \theta \in \mbox{CR}_5 \\
780.51 + 5.07\Delta w_1 + 1.93\Delta w_2, &  \theta \in \mbox{CR}_6 \\
\end{cases}
\end{equation}
Following the procedure in Section III-C, we can get the demand adjustment at equilibrium in each critical region as well. Take $\mbox{CR}_4$ as an example, we have $\forall [\Delta w_1,\Delta w_2]^\top \in \mbox{CR}_4$:
\begin{align}
\label{eq:approximate-strategy}
    \Delta d_1^* ~ &= 18.75 + 0.76 \Delta w_1 \nonumber\\
    \Delta d_2^* ~ & = -20.00 \nonumber\\
    \Delta d_3^* ~& = -3.75 + 0.24 \Delta w_1 + \Delta w_2
\end{align}

If we choose a scenario inside $\mbox{CR}_4$, say $[\Delta w_1,\Delta w_2]=[-10,-20]^\top$kW, then according to \eqref{eq:approximate-obj}, the approximate optimal cost is \$770.96 ($=837.26-2.01\times 10-2.31\times 20$); according to \eqref{eq:approximate-strategy}, the approximate demand adjustments at equilibrium are $\Delta d_1^*= 11.10$ kW, $\Delta d_2^*=-20$ kW, $\Delta d_3^* = -26.10$ kW. Both the optimal cost and demand adjustments are the same as the one we get by solving problem \eqref{eq:mp-LP-Matrix} directly with given parameters $\theta=[\Delta w_1,\Delta w_2]=[-10,-20]^\top$kW. This validates the effectiveness of our algorithm. The optimal cost of the original quadratic problem \eqref{eq:central} is \$767.24, so the relative error is only 0.48\%, showing the linearization is accurate enough. Moreover, the obtained $\Delta d_k^*, \forall k=1,2,3$ are the same as well. We calculate the flexibility requirement of each user, shown in TABLE \ref{tab:flex-req}. To accommodate the volatile renewable generation, all users need to have the ability to both increase and decrease their demand. User 3 requires the most flexibility since the value of $\overline{r}_3^d-\underline{r}_3^d$ is the largest.
\begin{table}[h]
        \renewcommand{\arraystretch}{1.3}
        \renewcommand{\tabcolsep}{1em}
        \centering
        \caption{Flexibility requirement of 5-bus system}
        \label{tab:flex-req}
        \begin{tabular}{cccc}
                \hline 
               Elastic demand & 1 & 2 & 3\\
               \hline
               $d_k$ (kW) & 230 & 170 & 200 \\
               $\underline{r}_k^d$ (kW) & -8.16 & -20.00 & -70.24\\
               $\overline{r}_k^d$ (kW) & 27.11 & 36.75 & 48.83 \\
             $\overline{r}_k^d-\underline{r}_k^d$ (kW) & 35.27 & 56.75 & 119.07\\
                \hline
        \end{tabular}
\end{table}

\begin{figure}[t]
\centering
\includegraphics[width=0.8\columnwidth]{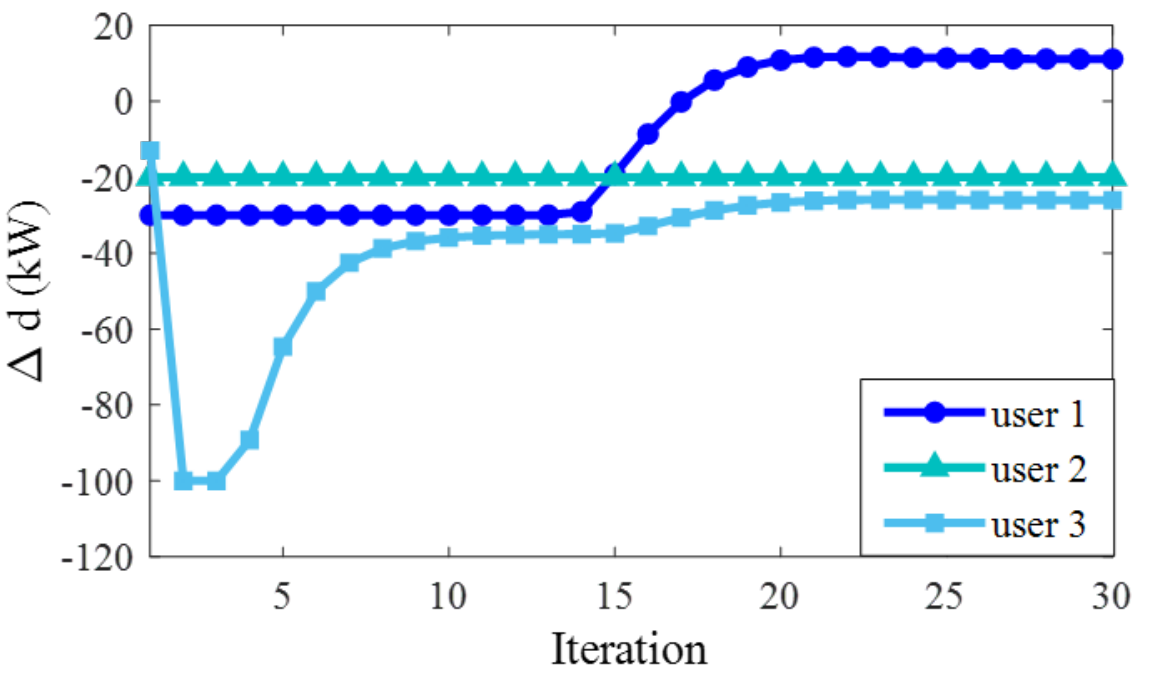}
\caption{Change of $\Delta d$ during iterations and the energy sharing market reaches an equilibrium after about 20 iterations.}
\label{fig:d}
\end{figure}

\begin{figure}[t]
\centering
\includegraphics[width=0.8\columnwidth]{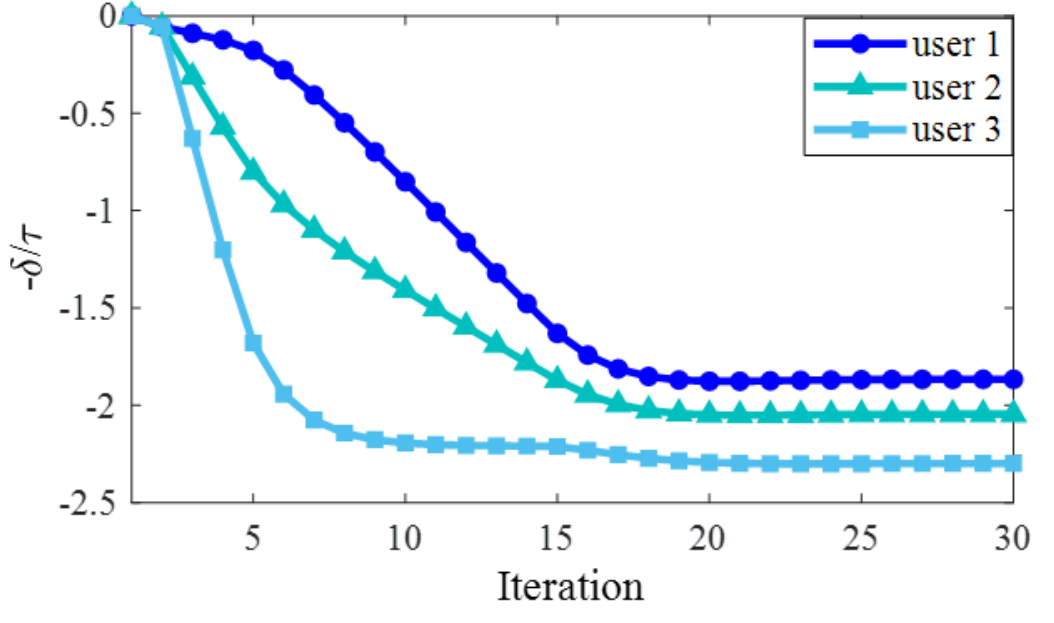}
\caption{Change of $\delta$ during iterations and $\frac{\delta^*}{\tau}=-\eta^*$ validates Proposition \ref{prop-1}.}
\label{fig:delta}
\end{figure}

Given $[\Delta w_1,\Delta w_2]=[-10,-20]^\top$kW, the equilibrium of the energy sharing game \eqref{eq:game1}-\eqref{eq:game2} can be reached via a best-response algorithm. The change of $\Delta d$ and $\delta$ during iterations are shown in Fig.\ref{fig:d} and Fig. \ref{fig:delta}, respectively. We can find that the demand adjustments converge to the above approximate values. The value of dual variables of problem \eqref{eq:central} at optimum can be retrieved by $\gamma_i$ in $\mbox{CR}_i$, i.e. $\eta_1^*=-1.92$, $\eta_2^*=-2.08$, $\eta_3^*=-2.31$, which are very close to the value of $-\delta/\tau$ plotted in Fig. \ref{fig:delta}. The above analysis validates Proposition \ref{prop-1}.


\subsection{69-bus system}
\begin{figure}[t]
\centering
\includegraphics[width=1.0\columnwidth]{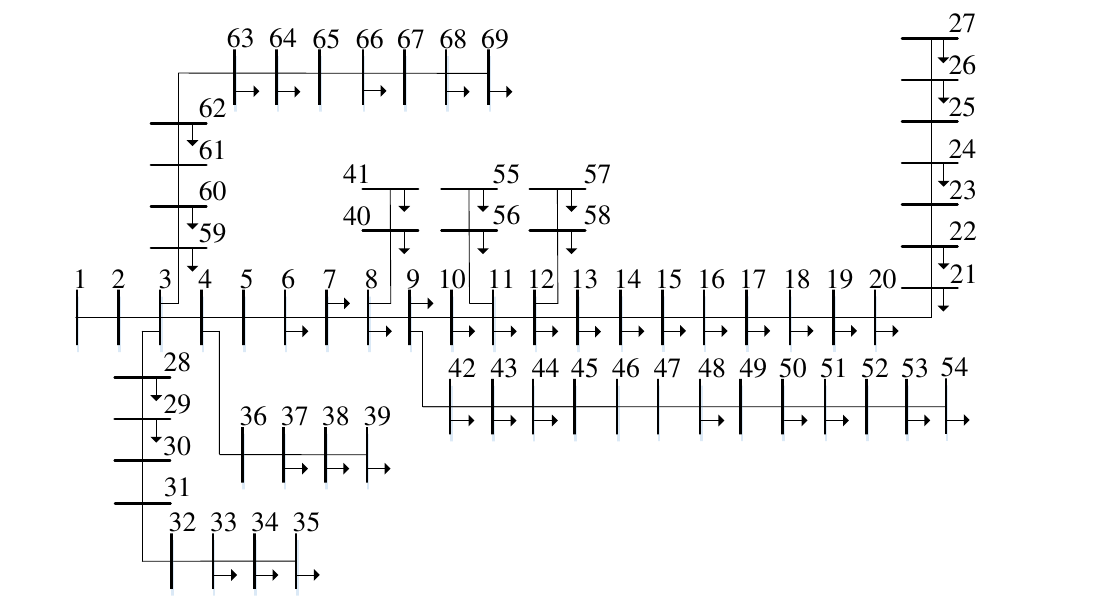}
\caption{Topology of the 69-bus system. Three wind farms are connected to nodes 9, 30, and 60; six elastic demands are located at nodes 12, 23, 32, 42, 53, and 62.}
\label{fig:69bus}
\end{figure}

\begin{figure}[t]
\centering
\includegraphics[width=0.75\columnwidth]{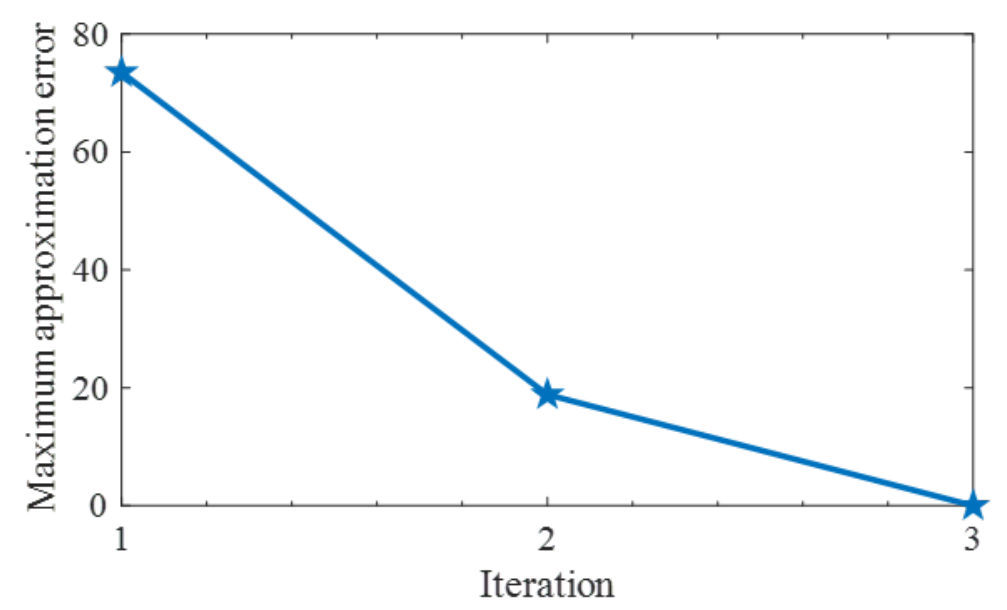}
\caption{The maximum approximation error of the AVG algorithm during iterations. The algorithm can output the critical regions after three iterations.}
\label{fig:error}
\end{figure}

\begin{figure}[t]
\centering
\includegraphics[width=0.75\columnwidth]{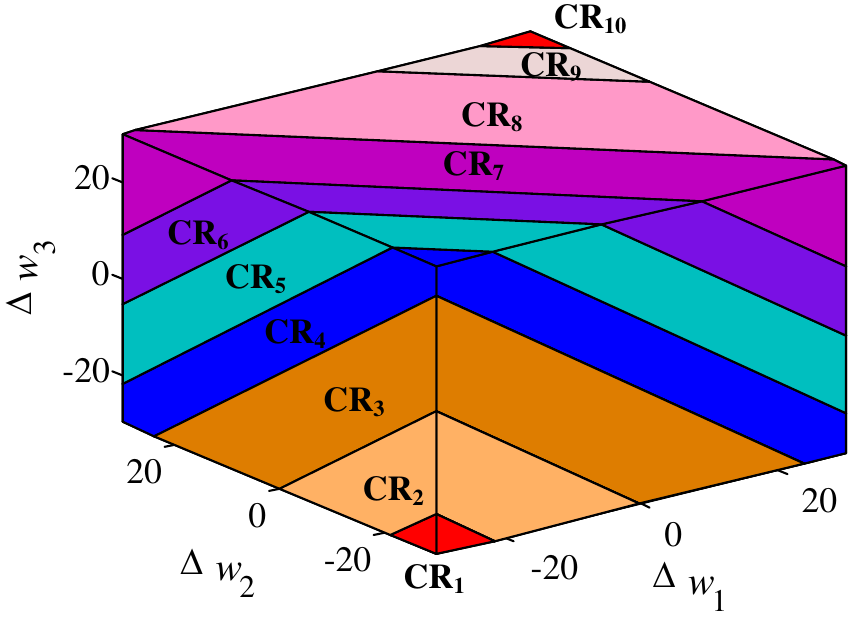}
\caption{Critical region of 69-bus system. The parameter range is divided into ten critical regions and within each region, the optimal cost is a linear function in $\theta$ given in \eqref{eq:approximate-obj2}.}
\label{fig:CR69}
\end{figure}

A larger case with a modified 69-bus microgrid \cite{kadir2013optimal} is tested with topology in Fig. \ref{fig:69bus}. There are three renewable generators connected to node 9, 30, and 60, respectively. Suppose their real-time output may deviate within $[-30,30]$ kW. Applying Algorithm 1, the change of the maximum approximation error (defined in \eqref{eq:Appr-Err-1}) during iterations is shown in Fig. \ref{fig:error}. The error decreases fast and the algorithm can output 10 critical regions (as in Fig. \ref{fig:CR69}) after three iterations. The approximate optimal cost function $\underline{v}(\theta)$ in each critical region is

\begin{equation}
\label{eq:approximate-obj2}
v(\theta) = \begin{cases}
353.05+[2.90, 2.90, 2.90]\theta
, &  \theta \in \mbox{CR}_1 \\
489.83+[4.58, 4.58, 4.58]\theta
, &  \theta \in \mbox{CR}_2 \\
527.73+[5.21, 5.21, 5.21]\theta
, &  \theta \in \mbox{CR}_3 \\
538.89+[5.52, 5.52, 5.52]\theta
, &  \theta \in \mbox{CR}_4 \\
549.33+[6.00, 6.00, 6.00]\theta
, &  \theta \in \mbox{CR}_5 \\
560.98+[8.09, 8.09, 8.09]\theta
, &  \theta \in \mbox{CR}_6 \\
548.24+[9.52, 9.52, 9.52]\theta
, &  \theta \in \mbox{CR}_7 \\
485.84+[11.47, 11.47, 11.47]\theta
, &  \theta \in \mbox{CR}_8 \\
361.11+[13.32, 13.32, 13.32]\theta
, &  \theta \in \mbox{CR}_9 \\
-9.46+[17.81, 17.81, 17.81]\theta
, &  \theta \in \mbox{CR}_{10} \\
\end{cases}
\end{equation}

To examine the accuracy of the approximation, let's take $\mbox{CR}_4$ as an example. Following the procedure in Section III-C, the demand adjustments at equilibrium are:
\begin{align}
\label{eq:approximate-strategy2}
    \Delta d_1^* = -17,~ \Delta d_2^*=-9,~\Delta d_3^*=-4,~\Delta d_4^*=8 \nonumber\\
    \Delta d_5^*=5, \Delta d_6^*=6+\Delta w_1+\Delta w_2+\Delta w_3
\end{align}

For scenario $\theta=[-30, -30, 30]$ kW $\in \mbox{CR}_4$, according to \eqref{eq:approximate-obj2} and \eqref{eq:approximate-strategy2}, the optimal cost is \$373.29 and $\Delta d_1^*=-17$ kW, $\Delta d_2^*=-9$ kW, $\Delta d_3^*=-4$ kW, $\Delta d_4^*=8$ kW, $\Delta d_5^*=5$ kW, $\Delta d_6^*=-24$ kW, which are the same as the one obtained by solving problem \eqref{eq:central}.
Similarly, we can calculate the flexibility requirement for each elastic demand as in TABLE \ref{tab:flex-req-69}. We notice that while elastic demands 2, 3, 4 still have some redundancy, the adjustable ability of elastic demands 1, 5, 6 has been fully exploited as $d_k+\underline{r}_k^d=\underline{D}_k$ and $d_k+\overline{r}_k^d=\overline{D}_k$, for $k=1, 5, 6$. Therefore, in long-term planning, we may put more emphasis on increasing the flexibility of demands 1, 5, 6.
\begin{table}[h]
        \renewcommand{\arraystretch}{1.3}
        \renewcommand{\tabcolsep}{1em}
        \centering
        \caption{Flexibility requirement of 69-bus system}
        \label{tab:flex-req-69}
        \begin{tabular}{ccccc}
                \hline 
                Elastic demand & $d_k$ (kW) & $\underline{r}_k^d$ (kW) & $\overline{r}_k^d$ (kW) & Range\\
               \hline
               1 & 40 & -30 & 10 & [10, 50]\\
               2 & 20 & -10.09 & 15 & [0, 35] \\
               3 & 30 & -25 & 24 & [5, 70]\\
               4 & 10 & -2 & 10 & [0, 20]\\
               5 & 10 & -5 &  10 & [5, 20]\\
               6 & 40 & -30 & 10 & [10, 50]\\
                \hline
        \end{tabular}
\end{table}

To show the applicability and scalability of the proposed algorithm, we test its performance under different numbers of renewable generators and users in TABLE \ref{tab:AVG-time}; and under different settings in TABLE \ref{tab:time}. The impact of the number of users on the computational time is marginal. The time needed increases with the number of RGs, but in this paper, we only consider the distributed renewable generators with a relatively large scale (say about 500kW) while the uncertainties from small-scale units are neglected. Usually, in a residential community, there are few such relatively large-scale units so our proposed method can still be applied. 
\begin{table}[h]
        \renewcommand{\arraystretch}{1.3}
        \renewcommand{\tabcolsep}{1em}
        \centering
        \caption{Computational time of AVG algorithm under different settings}
        \label{tab:AVG-time}
        \begin{tabular}{cccc}
                \hline 
               Cases & 3 RGs,69 users & 3 RGs,345 users &	3 RGs,690 users\\
               Time (s) & 103.28 & 122.34 & 92.37 \\
               \hline
               Cases & 6 RGs,69 users & 6 RGs,345 users	& 6 RGs,690 users \\
               Time (s) & 1507.94 &	1314.58 &	1644.23\\
                \hline
        \end{tabular}
\end{table}
\begin{table}[h]
        \renewcommand{\arraystretch}{1.3}
        \renewcommand{\tabcolsep}{1em}
        \centering
        \caption{Computational time}
        \label{tab:time}
        \begin{tabular}{cccc}
                \hline 
                Setting & 5-bus with $F_l$ & 5-bus with $2F_l$ & 5-bus with $4F_l$\\
               Time (s) & 22.09 & 13.89 & 10.31 \\
               \hline
               Setting & 69-bus with $F_l$ & 69-bus with $2F_l$ & 69-bus with $4F_l$ \\
               Time (s) & 103.28 & 60.32 & 64.19\\
                \hline
        \end{tabular}
\end{table}

To test the practicability of the proposed mechanism, we compare the computational time under different numbers of users in TABLE \ref{tab:AVG-user}. The computational time does not change much when the number of users increases since the decision-making for each user can be done in parallel. Our algorithm can handle a large number of users.
\begin{table}[h]
        \renewcommand{\arraystretch}{1.3}
        \renewcommand{\tabcolsep}{1em}
        \centering
        \caption{Computational time of the energy sharing mechanism}
        \label{tab:AVG-user}
        \begin{tabular}{cccc}
                \hline 
               Number of users & 69	& 345 & 690\\
               Time (s) & 29.74 & 30.06 & 30.21 \\
               \hline
        \end{tabular}
\end{table}

\section{Conclusion}
The prosperity of distributed renewable generators calls for effective management of end-users to hedge against uncertainty. In this paper, an innovative energy sharing mechanism is proposed that can facilitate the energy exchange among end-users with local information. The interaction among users and the operator is formulated as a generalized Nash game, whose equilibrium is proved to exist and can be obtained by a multi-parametric program parameterized in the renewable output deviations. To quantify the flexibility requirement of individual users, an alternative vertex generation algorithm is developed to output a piecewise linear estimator of both the optimal cost function and demand adjustment strategies at equilibrium via dual problem. It can avoid the impact of degeneracy since it does not require visiting the graph of degenerate base. Case studies validate the accuracy and scalability of our method.

Future research may further polish the proposed energy sharing mechanism by considering more realistic situations: 1) \emph{Bounded rationality of users.} In this paper, the users are assumed to be fully rational while in practice they may have bounded rationality. This can be analyzed with the prospect theory \cite{jhala2018prospect}. 2) \emph{AC power flow model.} Lossless DC power flow model is adopted in this paper for simplification, which shall be extended to incorporate AC power flow model with the help of convex relaxation \cite{huang2016sufficient} or linearization techniques \cite{chen2018energy}. 3) \emph{Diversity of devices.} Different types of flexible devices such as energy storage shall be integrated by changing the proposed model to a multiple time step version. More future research directions can be found in \cite{chen2021peer}.

\appendix
\makeatletter
\@addtoreset{equation}{section}
\@addtoreset{theorem}{section}
\makeatother
\setcounter{equation}{0}  
\renewcommand{\theequation}{A.\arabic{equation}}
\subsection{Proof of Proposition \ref{prop-1}}
\label{apen-1}  
Problem \eqref{eq:operator} can be equivalently written as
\begin{align}
\label{eq:apen-operator}
    \mathcal{\min}_{q_k^c,\forall k \in \mathcal{I} \cup \mathcal{J}}~ & \sum \nolimits_{k \in \mathcal{I} \cup \mathcal{J}}(q_k^c-q_k)^2 \nonumber\\
    \mbox{s.t.}~ &  q^c \in \mathcal{Q}
\end{align}
Suppose $(\Delta d^*,q^*,q^{c*}, \delta^*)$ is an GNE of the sharing game, then $q^{c*}$ is the optimal solution of \eqref{eq:apen-operator}. Therefore, we have
\begin{align}
\label{eq:condition-1}
    \sum \limits_{k \in \mathcal{I} \cup \mathcal{J}} (q_k^c-q_k^{c*})(q_k^{c*}-q_k^*) \ge 0, \forall q^c \in \mathcal{Q}
\end{align}
Also, $(\Delta d^*,q^*)$ is the optimal solution of problem \eqref{eq:consumer} or \eqref{eq:prosumer}, which is equivalent to 
\begin{align}
    \mathop{\min}_{\Delta d_k,b_k} ~ & \left\{\begin{aligned}
    & f_k(\Delta d_k) + \tau/2(q_k^{c*}+\delta_k^*-d_k-\Delta d_k)^2,\forall k \in \mathcal{I}\\
    & f_k(\Delta d_k) + \tau/2(q_k^{c*}+\delta_k^*-d_k-\Delta d_k+w_k+\Delta w_k)^2,\forall k \in \mathcal{J}
    \end{aligned} \right.\nonumber\\
    \mbox{s.t.} ~ & \Delta d_k \in \mathcal{D}_k
\end{align}
where $\mathcal{D}_k:=\{\Delta d_k ~|~ \underline{D}_k \le d_k+\Delta d_k \le \overline{D}_k\}$. Therefore,
\begin{align}
\label{eq:condition-2}
    f_k(\Delta d_k) - f_k(\Delta d_k^*) -\tau(\Delta d_k - \Delta d_k^*) \delta_k^* \ge 0, \forall k \in \mathcal{I}\cup \mathcal{J}
\end{align}
Similarly, let $(\Delta \hat d, \hat q^c)$ be the optimal solution of problem \eqref{eq:central} and $\hat \eta$ the value of the dual variable. Then, we have $\forall (\Delta d, q^c, \eta) \in \prod_{k \in \mathcal{I} \cup \mathcal{J}} \mathcal{D}_k \times \mathcal{Q} \times \mathbb{R}^{(I+J)}$:
\begin{align}
\label{eq:condition}
\left[f_k(\Delta d_k)- f_k(\Delta \hat d_k)+ (\Delta d_k-\Delta \hat d_k)\hat \eta_k\right] & ~\ge 0 ,\forall k \in \mathcal{I} \cup \mathcal{J} \nonumber\\
-\sum \limits_{k \in \mathcal{I} \cup \mathcal{J}} (q_k^c-\hat q_k^c)(\hat \eta_k) &~\ge 0 \nonumber\\
\sum \limits_{k \in \mathcal{I}}  (\eta_k-\hat \eta_k)(d_k+ \Delta \hat d_k-\hat q_k^c) &\nonumber\\
+\sum \limits_{k \in \mathcal{J}}  (\eta_k-\hat \eta_k)(d_k+ \Delta \hat d_k-w_k-\Delta w_k-\hat q_k^c) &~\le 0
\end{align}
Based on the above optimality conditions, we can analyze the existence and uniqueness of the GNE of the game $\mathcal{G}=\{\mathcal{K},\mathcal{X},\Pi\}$ as follows.

\emph{Existence}. When problem \eqref{eq:central} is feasible, since the objective function \eqref{eq:central.1} is strictly convex and the constraints constitute a close convex set, problem \eqref{eq:central} has a unique optimal solution $(\Delta \hat d, \hat q^c, \hat \eta)$. Let $\Delta d^*=\Delta \hat d$, $\delta^*=-\hat \eta/\tau$, $q^{c*}=\hat q^c$, $q^*=\hat q^c-\delta^*$. It is easy to check that condition \eqref{eq:condition-1}-\eqref{eq:condition-2} are satisfied so that $(\Delta d^*,q^*,q^{c*},\delta^*)$ is an GNE of the game $\mathcal{G}$.

\emph{Uniqueness}. Suppose $(\Delta d^*,q^*,q^{c*},\delta^*)$ is an GNE of the game $\mathcal{G}$, we have $q_k^{c*} = d_k+\Delta d^*_k,\forall k \in \mathcal{I}$, $q_k^{c*} = d_k+\Delta d_k^* - w_k - \Delta w_k,\forall k \in \mathcal{J}$ due to constraint \eqref{eq:consumer.2}, \eqref{eq:prosumer.2} and $q^{c*}=\delta^*+q^*$. Let $\Delta \hat d=\Delta d^*$, $\hat q^c=q^{c*}$, $\hat \eta= -\tau\delta^*$, then condition \eqref{eq:condition} is met. Therefore, $\Delta d^*$ is the optimal solution of \eqref{eq:central} and is unique.

\setcounter{equation}{0}
\renewcommand{\theequation}{B.\arabic{equation}}
\subsection{Redundancy Elimination}
Before introducing methods to eliminate redundancy, we first give the following definitions:
\label{apen-2}

\vspace{0.5em}
\noindent \textbf{Definition B1.} (Redundant Constraint)
A constraint is said to be \emph{redundant} if removing the constraint does not change the critical region. 
\vspace{0.5em}

\noindent \textbf{Definition B2.} (Minimal Representation)
Polyhedron (\ref{eq:Polyhedral-CR}) is said to be a minimal representation of CR$_i$ if all constraints are non-redundant.
\vspace{0.5em}

To remove redundant constraints to construct a minimal representation of CR$_i$, a method based on Nonhomogeneous Farkas Lemma is proposed.

\vspace{0.5em}
\noindent \textbf{Lemma B1.} (Nonhomogeneous Farkas Lemma \cite{xu1994generalization})
In the following two sets
\begin{equation*}
\begin{lgathered}
\mathcal{P}_1 = \{ u | A u \le t,~a^\top u > t^\prime \} \\
\mathcal{P}_2 = \{ v | A^\top v = a,~ t^\prime \ge v^\top t,~ v \ge 0 \} 
\end{lgathered}
\end{equation*}
Matrix multiplication is compatible in dimension. Then $\mathcal{P}_1$ is empty if and only if $\mathcal{P}_2$ is non-empty.
\vspace{0.5em}

Consider constraint $H_j \theta \le h_j$, if it is redundant, the polyhedron $\{\theta | H_{[-j]} \le h_{[-j]} \}$ defined by the remaining constraints must be a subset of $\{\theta | H_j \theta \le h_j\}$. In other words, $\{\theta | H_{[-j]} \le h_{[-j]} \} \cap \{\theta | H_j \theta > h_j\} = \emptyset$. So we have 

\vspace{0.5em}
\noindent \textbf{Theorem B1.}
Constraint $H_j \theta \le h_j$ in Polyhedron (\ref{eq:Polyhedral-CR}) is redundant if the following LP has a feasible solution
\begin{equation}
\label{eq:Redun-Test-LP}
v^\top  H_{[-j]} = H_j,~ h_j^\prime \ge v^\top h_{[-j]},~ v \ge 0
\end{equation}
\label{thm:Redundancy}

If all constraints in polyhedron (\ref{eq:Polyhedral-CR}) are screened, and redundant ones are removed, then the remaining constraints constitute a minimal representation of CR$_i$. Similarly, we can define and eliminate redundancy in $\underline v(\theta)$ as follows. 

\vspace{0.5em}
\noindent \textbf{Definition B3.} (Redundant Piece)
 In the piecewise linear function (\ref{eq:OVF-Dual-L}), a piece $\gamma^{\top}_i b+ \gamma^{\top}_i B \theta$ is said to be redundant if it never reaches maximum for all $\theta \in \Theta$.
\vspace{0.5em}

\noindent {\bf Claim B1.}  The $j$-th piece  in $\underline v(\theta)$\ is redundant if the $j$-th constraint in the polyhedron $\mbox{epi}[\underline v(\theta)]$ is redundant. The epigraph of the piecewise linear function (\ref{eq:OVF-Dual-L}) is the following polyhedron 
\begin{equation}
\label{eq:Epi-OVF}
\mbox{epi}[\underline v(\theta)] = \left\{ (\theta,\kappa) \in \mathbb R^{p+1} \middle| \begin{gathered}
\kappa \ge \gamma^{\top}_1 t+ \gamma^{\top}_1 B \theta  \\
\vdots  \\
\kappa \ge \gamma^{\top}_n t+ \gamma^{\top}_n B \theta 
\end{gathered} \right\}
\end{equation}
\vspace{0.5em}

Similarly, according to Theorem B1, the redundancy in $\underline v(\theta)$ can be removed.

\ifCLASSOPTIONcaptionsoff
\newpage
\fi
        
        \bibliographystyle{IEEEtran}
        \bibliography{IEEEabrv,mybib}

\begin{thebibliography}{10}
\providecommand{\url}[1]{#1}
\csname url@samestyle\endcsname
\providecommand{\newblock}{\relax}
\providecommand{\bibinfo}[2]{#2}
\providecommand{\BIBentrySTDinterwordspacing}{\spaceskip=0pt\relax}
\providecommand{\BIBentryALTinterwordstretchfactor}{4}
\providecommand{\BIBentryALTinterwordspacing}{\spaceskip=\fontdimen2\font plus
\BIBentryALTinterwordstretchfactor\fontdimen3\font minus
  \fontdimen4\font\relax}
\providecommand{\BIBforeignlanguage}[2]{{%
\expandafter\ifx\csname l@#1\endcsname\relax
\typeout{** WARNING: IEEEtran.bst: No hyphenation pattern has been}%
\typeout{** loaded for the language `#1'. Using the pattern for}%
\typeout{** the default language instead.}%
\else
\language=\csname l@#1\endcsname
\fi
#2}}
\providecommand{\BIBdecl}{\relax}
\BIBdecl

\bibitem{chen2020decentralized}
Y.~Chen, T.~Li, C.~Zhao, and W.~Wei, ``Decentralized provision of renewable
  predictions within a virtual power plant,'' \emph{IEEE Transactions on Power
  Systems}, vol.~36, no.~3, pp. 2652--2662, 2021.

\bibitem{shuai2018stochastic}
H.~Shuai, J.~Fang, X.~Ai, Y.~Tang, J.~Wen, and H.~He, ``Stochastic optimization
  of economic dispatch for microgrid based on approximate dynamic
  programming,'' \emph{IEEE Transactions on Smart Grid}, vol.~10, no.~3, pp.
  2440--2452, 2018.

\bibitem{jiang2011robust}
R.~Jiang, J.~Wang, and Y.~Guan, ``Robust unit commitment with wind power and
  pumped storage hydro,'' \emph{IEEE Transactions on Power Systems}, vol.~27,
  no.~2, pp. 800--810, 2011.

\bibitem{wei2015distributionally}
W.~Wei, F.~Liu, and S.~Mei, ``Distributionally robust co-optimization of energy
  and reserve dispatch,'' \emph{IEEE Transactions on Sustainable Energy},
  vol.~7, no.~1, pp. 289--300, 2015.

\bibitem{bertsimas2012adaptive}
D.~Bertsimas, E.~Litvinov, X.~A. Sun, J.~Zhao, and T.~Zheng, ``Adaptive robust
  optimization for the security constrained unit commitment problem,''
  \emph{IEEE Transactions on Power Systems}, vol.~28, no.~1, pp. 52--63, 2012.

\bibitem{zeng2013solving}
B.~Zeng and L.~Zhao, ``Solving two-stage robust optimization problems using a
  column-and-constraint generation method,'' \emph{Operations Research
  Letters}, vol.~41, no.~5, pp. 457--461, 2013.

\bibitem{chen2018analyzing}
Y.~Chen, W.~Wei, F.~Liu, Q.~Wu, and S.~Mei, ``Analyzing and validating the
  economic efficiency of managing a cluster of energy hubs in multi-carrier
  energy systems,'' \emph{Applied energy}, vol. 230, pp. 403--416, 2018.

\bibitem{chen2021peer}
Y.~Chen and C.~Zhao, ``Peer-to-peer energy sharing: A new business model
  towards a low-carbon future,'' \emph{arXiv preprint arXiv:2108.04057}, 2021.

\bibitem{tushar2020coalition}
W.~Tushar, T.~K. Saha, C.~Yuen, M.~I. Azim, T.~Morstyn, H.~V. Poor, D.~Niyato,
  and R.~Bean, ``A coalition formation game framework for peer-to-peer energy
  trading,'' \emph{Applied Energy}, vol. 261, p. 114436, 2020.

\bibitem{azim2021coalition}
M.~I. Azim, W.~Tushar, and T.~K. Saha, ``Coalition graph game-based p2p energy
  trading with local voltage management,'' \emph{IEEE Transactions on Smart
  Grid}, 2021.

\bibitem{tushar2019grid}
W.~Tushar, T.~K. Saha, C.~Yuen, T.~Morstyn, H.~V. Poor, R.~Bean \emph{et~al.},
  ``Grid influenced peer-to-peer energy trading,'' \emph{IEEE Transactions on
  Smart Grid}, vol.~11, no.~2, pp. 1407--1418, 2019.

\bibitem{liu2017energy}
N.~Liu, X.~Yu, C.~Wang, and J.~Wang, ``Energy sharing management for microgrids
  with {PV} prosumers: A stackelberg game approach,'' \emph{IEEE Transactions
  on Industrial Informatics}, vol.~13, no.~3, pp. 1088--1098, 2017.

\bibitem{tushar2018peer}
W.~Tushar, T.~K. Saha, C.~Yuen, P.~Liddell, R.~Bean, and H.~V. Poor,
  ``Peer-to-peer energy trading with sustainable user participation: A game
  theoretic approach,'' \emph{IEEE Access}, vol.~6, pp. 62\,932--62\,943, 2018.

\bibitem{liu2018energy}
N.~Liu, M.~Cheng, X.~Yu, J.~Zhong, and J.~Lei, ``Energy-sharing provider for
  {PV} prosumer clusters: A hybrid approach using stochastic programming and
  stackelberg game,'' \emph{IEEE Transactions on Industrial Electronics},
  vol.~65, no.~8, pp. 6740--6750, 2018.

\bibitem{xu2020two}
X.~Xu, J.~Li, Y.~Xu, Z.~Xu, and C.~S. Lai, ``A two-stage game-theoretic method
  for residential {PV} panels planning considering energy sharing mechanism,''
  \emph{IEEE Transactions on Power Systems}, vol.~35, no.~5, pp. 3562--3573,
  2020.

\bibitem{li2017risk}
C.~Li, Y.~Xu, X.~Yu, C.~Ryan, and T.~Huang, ``Risk-averse energy trading in
  multienergy microgrids: A two-stage stochastic game approach,'' \emph{IEEE
  Transactions on Industrial Informatics}, vol.~13, no.~5, pp. 2620--2630,
  2017.

\bibitem{chen2019energy}
Y.~Chen, S.~Mei, F.~Zhou, S.~H. Low, W.~Wei, and F.~Liu, ``An energy sharing
  game with generalized demand bidding: Model and properties,'' \emph{IEEE
  Transactions on Smart Grid}, vol.~11, no.~3, pp. 2055--2066, 2019.

\bibitem{chen2020approaching}
Y.~Chen, C.~Zhao, S.~H. Low, and S.~Mei, ``Approaching prosumer social optimum
  via energy sharing with proof of convergence,'' \emph{IEEE Transactions on
  Smart Grid}, vol.~12, no.~3, pp. 2484--2495, 2021.

\bibitem{cui2019peer}
S.~Cui, Y.-W. Wang, and J.-W. Xiao, ``Peer-to-peer energy sharing among smart
  energy buildings by distributed transaction,'' \emph{IEEE Transactions on
  Smart Grid}, vol.~10, no.~6, pp. 6491--6501, 2019.

\bibitem{cui2019efficient}
S.~Cui, Y.-W. Wang, Y.~Shi, and J.-W. Xiao, ``An efficient peer-to-peer
  energy-sharing framework for numerous community prosumers,'' \emph{IEEE
  Transactions on Industrial Informatics}, vol.~16, no.~12, pp. 7402--7412,
  2019.

\bibitem{zhao2014variable}
J.~Zhao, T.~Zheng, and E.~Litvinov, ``Variable resource dispatch through
  do-not-exceed limit,'' \emph{IEEE Transactions on Power Systems}, vol.~30,
  no.~2, pp. 820--828, 2014.

\bibitem{wei2014dispatchable}
W.~Wei, F.~Liu, and S.~Mei, ``Dispatchable region of the variable wind
  generation,'' \emph{IEEE Transactions on Power Systems}, vol.~30, no.~5, pp.
  2755--2765, 2014.

\bibitem{li2016data}
Z.~Li, F.~Qiu, and J.~Wang, ``Data-driven real-time power dispatch for
  maximizing variable renewable generation,'' \emph{Applied Energy}, vol. 170,
  pp. 304--313, 2016.

\bibitem{ma2019distributionally}
H.~Ma, R.~Jiang, and Z.~Yan, ``Distributionally robust co-optimization of power
  dispatch and do-not-exceed limits,'' \emph{IEEE Transactions on Power
  Systems}, vol.~35, no.~2, pp. 887--897, 2019.

\bibitem{korad2015enhancement}
A.~S. Korad and K.~W. Hedman, ``Enhancement of do-not-exceed limits with robust
  corrective topology control,'' \emph{IEEE Transactions on Power Systems},
  vol.~31, no.~3, pp. 1889--1899, 2015.

\bibitem{nosair2015reconstructing}
H.~Nosair and F.~Bouffard, ``Reconstructing operating reserve: Flexibility for
  sustainable power systems,'' \emph{IEEE Transactions on Sustainable Energy},
  vol.~6, no.~4, pp. 1624--1637, 2015.

\bibitem{dvorkin2014assessing}
Y.~Dvorkin, D.~S. Kirschen, and M.~A. Ortega-Vazquez, ``Assessing flexibility
  requirements in power systems,'' \emph{IET Generation, Transmission \&
  Distribution}, vol.~8, no.~11, pp. 1820--1830, 2014.

\bibitem{wei2020nodal}
W.~Wei, Y.~Chen, C.~Wang, Q.~Wu, and M.~Shahidehpour, ``Nodal flexibility
  requirements for tackling renewable power fluctuations,'' \emph{IEEE
  Transactions on Power Systems}, 2020.

\bibitem{chen2021energy}
Y.~Chen, W.~Wei, H.~Wang, Q.~Zhou, and J.~P. Catal{\~a}o, ``An energy sharing
  mechanism achieving the same flexibility as centralized dispatch,''
  \emph{IEEE Transactions on Smart Grid}, 2021.

\bibitem{filippi2004algorithm}
C.~Filippi, ``An algorithm for approximate multiparametric linear
  programming,'' \emph{Journal of optimization theory and applications}, vol.
  120, no.~1, pp. 73--95, 2004.

\bibitem{bemporad2006algorithm}
A.~Bemporad and C.~Filippi, ``An algorithm for approximate multiparametric
  convex programming,'' \emph{Computational optimization and applications},
  vol.~35, no.~1, pp. 87--108, 2006.

\bibitem{ghorani2019main}
R.~Ghorani, M.~Fotuhi-Firuzabad, and M.~Moeini-Aghtaie, ``Main challenges of
  implementing penalty mechanisms in transactive electricity markets,''
  \emph{IEEE Transactions on Power Systems}, vol.~34, no.~5, pp. 3954--3956,
  2019.

\bibitem{pei2016optimal}
W.~Pei, Y.~Du, W.~Deng, K.~Sheng, H.~Xiao, and H.~Qu, ``Optimal bidding
  strategy and intramarket mechanism of microgrid aggregator in real-time
  balancing market,'' \emph{IEEE Transactions on Industrial Informatics},
  vol.~12, no.~2, pp. 587--596, 2016.

\bibitem{wei2014robust}
W.~Wei, F.~Liu, S.~Mei, and Y.~Hou, ``Robust energy and reserve dispatch under
  variable renewable generation,'' \emph{IEEE Transactions on Smart Grid},
  vol.~6, no.~1, pp. 369--380, 2014.

\bibitem{nimbargi2016review}
S.~Nimbargi, S.~Mhaisne, S.~Nangare, and M.~Sinha, ``Review on ami technology
  for smart meter,'' in \emph{2016 IEEE International Conference on Advances in
  Electronics, Communication and Computer Technology (ICAECCT)}.\hskip 1em plus
  0.5em minus 0.4em\relax IEEE, 2016, pp. 21--27.

\bibitem{harker1991generalized}
P.~T. Harker, ``Generalized nash games and quasi-variational inequalities,''
  \emph{European journal of Operational research}, vol.~54, no.~1, pp. 81--94,
  1991.

\bibitem{facchinei2010generalized}
F.~Facchinei and C.~Kanzow, ``Generalized nash equilibrium problems,''
  \emph{Annals of Operations Research}, vol. 175, no.~1, pp. 177--211, 2010.

\bibitem{wu2011tighter}
L.~Wu, ``A tighter piecewise linear approximation of quadratic cost curves for
  unit commitment problems,'' \emph{IEEE Transactions on Power Systems},
  vol.~26, no.~4, pp. 2581--2583, 2011.

\bibitem{borrelli2003geometric}
F.~Borrelli, A.~Bemporad, and M.~Morari, ``Geometric algorithm for
  multiparametric linear programming,'' \emph{Journal of optimization theory
  and applications}, vol. 118, no.~3, pp. 515--540, 2003.

\bibitem{fortuny1981representation}
J.~Fortuny-Amat and B.~McCarl, ``A representation and economic interpretation
  of a two-level programming problem,'' \emph{Journal of the operational
  Research Society}, vol.~32, no.~9, pp. 783--792, 1981.

\bibitem{Data}
\emph{[Online] Available:},
  \url{https://sites.google.com/site/yuechenthu/data-sheet}.

\bibitem{kadir2013optimal}
A.~F.~A. KADIR, A.~Mohamed, H.~Shareef, and M.~Z.~C. WANIK, ``Optimal placement
  and sizing of distributed generations in distribution systems for minimizing
  losses and {THD}\_v using evolutionary programming,'' \emph{Turkish Journal
  of Electrical Engineering \& Computer Sciences}, vol.~21, no. Sup. 2, pp.
  2269--2282, 2013.

\bibitem{jhala2018prospect}
K.~Jhala, B.~Natarajan, and A.~Pahwa, ``Prospect theory-based active consumer
  behavior under variable electricity pricing,'' \emph{IEEE Transactions on
  Smart Grid}, vol.~10, no.~3, pp. 2809--2819, 2018.

\bibitem{huang2016sufficient}
S.~Huang, Q.~Wu, J.~Wang, and H.~Zhao, ``A sufficient condition on convex
  relaxation of ac optimal power flow in distribution networks,'' \emph{IEEE
  Transactions on Power Systems}, vol.~32, no.~2, pp. 1359--1368, 2016.

\bibitem{chen2018energy}
Y.~Chen, W.~Wei, F.~Liu, E.~E. Sauma, and S.~Mei, ``Energy trading and market
  equilibrium in integrated heat-power distribution systems,'' \emph{IEEE
  Transactions on Smart Grid}, vol.~10, no.~4, pp. 4080--4094, 2018.

\bibitem{xu1994generalization}
Z.~Xu, ``Generalization of nonhomogeneous farkas' lemma and applications,''
  \emph{Journal of Mathematical Analysis and Applications}, vol. 186, no.~3,
  pp. 726--734, 1994.

\end{thebibliography}

\end{document}